\documentclass{commat}

\title{Noncommutative algebra and representation theory: symmetry, structure  invariants}

\author{Samuel A. Lopes}

\authorinfo[%
S. Lopes]{
CMUP, Departamento de Matem\'atica, Faculdade de Ci\^encias, Universidade do Porto, Rua do Campo Alegre s/n, 4169--007 Porto, Portugal.}{%
slopes@fc.up.pt
}

\abstract{%
This is an abridged version of our Habilitation thesis. In these notes, we aim to summarize our research interests and achievements as well as motivate what drives our work: symmetry, structure and invariants. The paradigmatic example which permeates and often inspires our research is the Weyl algebra $\wa$. 
}



\usepackage{amssymb}
\usepackage{amsfonts}
\usepackage{amsmath}
\usepackage{amsthm}
\usepackage[dvipsnames]{xcolor}
\usepackage[shortlabels]{enumitem}
\usepackage{mathtools}
\usepackage{hyperref}
\usepackage[all]{xy}

\usepackage{accents}

\usepackage{tikz-cd}
\usepackage{faktor}

\usepackage[vcentermath]{youngtab}
\usepackage{young}
\usepackage{pgfplots}
\pgfplotsset{compat=1.18}

\usepackage{multicol}
\usepackage{array}
\usepackage{booktabs}

\newcommand\q{\quad}
\newcommand\ds{\displaystyle}

\newcommand{\overbar}[1]{\mkern 1mu\overline{\mkern-1mu#1\mkern-1mu}\mkern 1mu}
\newcommand{\ann}{\hbox{\rm Ann}}
\newcommand\tr{\mathsf{tr}}

\newcommand{\FF}{\mathbb{F}}  
\newcommand{\wa}[1][1]{\mathbb{A}_{#1}}  
\newcommand{\ZZ}{\mathbb{Z}}  
\newcommand{\CC}{\mathbb{C}}  
\newcommand{\KK}{\mathbb{K}} 

\newcommand{\A}{\mathsf{A}} 
\newcommand{\D}{\mathsf{D}}
\newcommand{\W}{\mathsf{W}} 

\newcommand{\pr}{\mathsf{u}}

\newcommand\ot{\otimes}

\newcommand\End{\mathsf {End}}
\newcommand\Hom{\mathsf{Hom}}

\newcommand\chara{\mathsf {char}}
\newcommand{\barr}[1]{\overline{#1}}
\newcommand{\eval}[2]{\left. #1 \right|_{#2}}

\newcommand{\hqfg}[1][f,g]{\mathcal{H}_q(#1)}
\newcommand{\hqqfg}[1][f',g']{\mathcal{H}_{q'}(#1)}
\newcommand{\hh}{\mathcal{H}}
\newcommand{\Sf}[1][f]{\mathsf{S}_{#1}}
\newcommand{\Tf}[1][q, g, \lambda]{\mathsf{T}_{#1}}
\newcommand{\Af}[1][\lambda, \mu]{\mathsf{A}_{q, f,g}(#1)}
\newcommand{\Bf}[1][\lambda, \mu]{\mathsf{B}_{q, f,g}(#1)}
\newcommand{\Cf}[1][\alpha]{\mathsf{C}_{q, f,g}(#1)}

\newcommand\inder{ \mathsf {Inder_\FF}}
\newcommand\der[1][\FF]{ \mathsf{Der}_{#1}}
\newcommand\aut[1][\FF]{ \mathsf{Aut}_{#1}}
\newcommand\hoch{\mathsf{HH}}
\newcommand\im{\mathsf{im}}
\newcommand\modd[1]{\ (\mathsf{mod} \, #1)}
\newcommand\degg{\, \mathsf{deg} \,}
\newcommand\dimm{\, \mathsf{dim}\,}
\newcommand\spann{\, \mathsf{span}}
\newcommand{\g}{\mathfrak g}
\newcommand{\h}{\mathfrak h}
\newcommand{\fsl}{\mathfrak{sl}}

\newcommand\ad{\mathsf{ad}}

\newcommand{\lb}[1]{\left[ #1\right]}
\newcommand{\pb}[1]{\left\{ #1\right\}}
\newcommand{\seq}[1]{\left( #1\right)}
\newcommand{\ip}[1]{\left\langle #1\right\rangle}

\usepackage{stmaryrd} 
\newcommand{\formal}[1]{\left\llbracket #1\right\rrbracket}

\definecolor{azulDM}{RGB}{130,167,199}
\definecolor{vermDM}{RGB}{203,88,56}
\newcommand{\red}{\bf \color{vermDM}} 
\newcommand\cred[1]{{\color{vermDM}#1}}

\VOLUME{32}
\NUMBER{3}
\YEAR{2024}
\firstpage{63}
\DOI{https://doi.org/10.46298/cm.11678}

\begin{document}

{
	\tableofcontents
	
	\thispagestyle{empty}
	
	\clearpage
}

\bigskip

\begin{center}
	\scalebox{.8}{$$\xymatrix{&(4,1)_{\red 11} \ar@{-}[dl]^4 \ar@{-}[dr]^1 &(3,2)_{\red 15} \ar@{-}[d]^1 \ar@{-}[dr]^2& (3,1^2)_{\red 32} \ar@{-}[dl]^3 \ar@{-}[dr]^1  &(2,2,1)_{\red 34}   & (2,1^3)_{\red 26} \ar@{-}[dl]^2 \ar@{-}[dr]^4 \\ 
			(4)_{\red 1} \ar@{-}[dr]^1 && (3,1)_{\red 7} \ar@{-}[dl]^3 \ar@{-}[dr]^1 & (2,2)_{\red 4} \ar@{-}[ur]^3 \ar@{-}[d]^1 & (2,1^2)_{\red 11} \ar@{-}[u]_2 \ar@{-}[dl]^2 \ar@{-}[dr]^3 && (1^4)_{\red 1} \ar@{-}[dl]_1  \\ 
			& (3)_{\red 1} \ar@{-}[dr]^1 && (2,1)_{\red 4} \ar@{-}[dl]^2 \ar@{-}[dr]_2 && (1^3)_{\red 1} \ar@{-}[dl]_1 \\ 
			&& (2)_{\red 1} \ar@{-}[dr]^1  && (1^2)_{\red 1} \ar@{-}[dl]_1 \\
			&&& (1)_{\red 1}}$$} 
\end{center}

\section{Introduction}\label{S:intro}

\begin{flushright}{\it
		``One can see the world with the $p$-eye and one can look at it with the $q$-eye. But if one wants to open both eyes at the same time, one goes crazy.''} \\
	Letter from Wolfgang Pauli to Werner Heisenberg, 1926
\end{flushright}
(Above, $p$ refers to the momentum and $q$ to the position of a particle in quantum mechanics and Pauli is alluding to Heisenberg's uncertainty principle.)

\bigskip

In these notes we aim to summarize our research interests and achievements as well as motivate what drives our work: symmetry, structure and invariants. The paradigmatic example which permeates and often inspires our research is the Weyl algebra $\wa$. When defined over the field $\CC$ of complex numbers, $\wa$ is the algebra of differential operators in one variable with polynomial coefficients. In other words, it is the associative unital subalgebra of the algebra $\End_\CC (\CC[x])$ of linear operators on the polynomial algebra $\CC[x]$ generated by $\CC[x]$ (viewed as left multiplication operators) and by its Lie algebra of derivations $\der[\CC](\CC[x])$. More generally, the $n$-th Weyl algebra over $\CC$, denoted by $\wa[n]$, is the algebra of differential operators in $n$ variables with polynomial coefficients. One can see that
\begin{equation}\label{E:nweyltensor}
	\wa[n]\simeq\underbrace{\wa\ot\cdots\ot\wa}_{n}
\end{equation}
so, in a sense, $\wa$ is the basic building block of $\wa[n]$.

The Weyl algebras are central in the theory of $\mathcal D$-modules (developed by Sato, Kashiwara and Bernstein in connection with the theory of linear partial differential equations), in quantum Physics, in algebraic geometry, linking commutative and noncommutative algebra, and of course in representation theory and the theory of noncommutative rings and algebras. 

By work of Dixmier (see e.g.\ \cite[Thm.\ 4.7.9]{jD96}), the Weyl algebras control the representation theory of the finite-dimensional complex nilpotent Lie algebras but they are also related to the representations of semisimple Lie algebras. Weyl algebras are thus of paramount importance in Lie theory and representation theory. They are also a prototypical example of a simple noncommutative Noetherian ring which is not a matrix ring over a division algebra.

Next we recall results and open problems about the Weyl algebra which will motivate and serve as guidelines for the lines of research expounded in this lecture.

\subsection{Conventions and notation}\label{SS:intro:cn}

Unless otherwise noted, all vector spaces and algebras are considered over an arbitrary field $\FF$, with algebraic closure $\barr \FF$ and group of units $\FF^*$. All rings and algebras are assumed to be unital and associative. To avoid deciding whether $0$ is a natural number or not (although it seems to us that it should be!), we use $\ZZ_{\geq 0}$ and $\ZZ_{> 0}$ to denote the sets of nonnegative and positive integers, respectively.

Given an associative algebra $A$ and elements $a,b \in A$, we use the commutator notation $[a,b]=ab-ba$ and let $\ad_a$ denote the endomorphism of $A$ defined by $\ad_a(b)=[a,b]$. The center of $A$ and the centralizer of an element $a\in A$ will be denoted by $\mathsf{Z}(A)$ and $\mathsf{C}_{A}(a)$, respectively. An element $c\in A$ is \textit{normal} if $cA=Ac$ (an ideal of $A$). We remark that the set of normal elements of $A$ forms a multiplicative monoid. If $\g$ is a Lie algebra, then we denote its universal enveloping algebra by $U(\g)$.

Unadorned $\ot$ will always mean $\ot_\FF$. For any set $E$, $1_E$ will denote the identity map on $E$.
Given  $f\in\FF[x]$, 
$f^{(k)}$ stands for the $k$-th derivative of $f$ with respect to $x$, which we also denote by $f'$ and $f''$ in case $k=1, 2$, respectively. If $f, g\in\FF[x]$ are not both zero, then we tacitly assume that $\gcd(f, g)$ is monic.

\subsection{The Weyl algebra: definition}\label{SS:intro:wa:def}

The (first) Weyl algebra $\wa(\FF)$ is the unital associative $\FF$-algebra with generators $x$ and $y$, subject to the relation $yx-xy=1$.
In other words, 
\begin{align*}
	\wa(\FF)=\FF\langle x, y\rangle/\seq{yx-xy-1},
\end{align*}
the quotient of the free unital associative algebra on generators $x$ and $y$, by the two-sided ideal generated by $yx-xy-1$.

In case $\chara(\FF)=0$, then (up to isomorphism),
\begin{align}\label{E:intro:wa:def:1}
	\wa(\FF)=\spann_\FF\pb{t^i\frac{d^j}{dt^j}\mid i, j\geq 0},
\end{align}
where the right-hand side is seen as linear operators on $\FF[t]$, under composition. In case $\chara(\FF)=p>0$, then $\frac{d^j}{dt^j}=0$ for all $j\geq p$, as an operator on $\FF[t]$, and $y^p$ is central in $\wa(\FF)$; we have, in this case, 
\begin{align}\label{E:intro:wa:def:2}
	\wa(\FF)/\seq{y^p}=\spann_\FF\pb{t^i\frac{d^j}{dt^j}\mid i, j\geq 0}.
\end{align}

\subsection{The Weyl algebra: representations}\label{SS:intro:wa:reps}

Throughout this subsection, assume that $\FF$ is algebraically closed.

\subsubsection{\texorpdfstring{$\chara(\FF)=0$}{char(F)=0}}

If $\chara(\FF)=0$, then the relation $yx-xy=1$ implies that $\wa(\FF)$ has no nonzero finite-dimen\-sional representations. Indeed, if $X$ and $Y$ are linear operators on a finite-dimensional vector space $V$, then 
\begin{align*}
	\tr\seq{YX-XY}=0\quad\text{whereas}\quad \tr (1_V)=\dim_\FF V.
\end{align*} 
However, \eqref{E:intro:wa:def:1} gives faithful simple representation of $\wa(\FF)$ on the polynomial algebra $\FF[t]$, where $x$ acts by multiplication by $t$ and $y$ as $\frac{d}{dt}$.

\begin{example}[Representations of the Heisenberg Lie algebra]
	Let $\mathfrak h$ be the $3$-dimensional Heisenberg Lie algebra. Then $\mathfrak h$ is spanned by elements $x, y, z$ with Lie brackets
	\begin{align*}
		[y, x]=z, \quad \lb{z,\mathfrak h}=0.
	\end{align*}
	The universal enveloping algebra of $\mathfrak h$, is the unital associative algebra $U(\mathfrak h)$ generated by $x, y, z$, with relations
	\begin{align*}
		yx-xy=z,\quad zx=xz, \quad zy=yz.
	\end{align*}
	Assume that $\FF$ is algebraically closed and of characteristic $0$. Then the simple finite-dimensional $\mathfrak h$-modules are $1$-dimensional (by Lie's theorem). This implies that if $V$ is a simple infinite-dimensional $\mathfrak h$-module, equivalently a $U(\mathfrak h)$-module, then $z$ will act on $V$ by some nonzero scalar $\alpha\in\FF^*$, so $V$ can be seen equivalently as an $U(\mathfrak h)/(z-\alpha)$-module. But $U(\mathfrak h)/(z-\alpha)\simeq \wa(\FF)$, so there is a bijective correspondence between simple infinite-dimensional representations of $\h$ and simple representations of $\wa(\FF)$.
\end{example}

\begin{flushright}{\it
		``But a deeper study reveals the existence of an enormous number of irreducible representations of $\h$, even for $\alpha$ ($\neq 0$) fixed. It seems that these representations defy classification. A similar phenomenon exists for $\g=\fsl_2$, and most certainly for all non-commutative Lie algebras.''} \\
	J.\ Dixmier, \cite[Preface]{jD77book}
\end{flushright}

In \cite{rB81}, Block undertook a remarkable and comprehensive study of the irreducible modules for the Weyl algebra $\wa(\FF)$ (and hence for the $3$-dimensional Heisenberg Lie algebra) and for the universal enveloping algebras of 
$\mathfrak{sl}_2$ and of the two-dimensional solvable Lie algebra over $\FF$. (Compare also \cite{AP74} for the $\mathfrak{sl}_2$ case.) Block's results were extended by Bavula in \cite{vB99} to more general Ore extensions over Dedekind domains, and by Bavula and van Oystaeyen in  \cite{BvO97}  to develop a representation theory for generalized Weyl algebras over Dedekind domains.  

\subsubsection{\texorpdfstring{$\chara(\FF)=p>0$}{char(F)=p>0}}

In this case, $\wa(\FF)$ is a free module over its center $\FF[x^p, y^p]$ with rank $p^2$ and it follows that all simple $\wa(\FF)$-modules are finite dimensional of dimension $p$. One such simple module is
\begin{align*}
	\FF[t]/t^p\FF[t],
\end{align*}
with the action induced from the natural action on $\FF[t]$.

\subsection{The Weyl algebra: automorphisms and derivations}\label{SS:intro:wa:ad}

Let $\g$ be a Lie algebra over $\FF$ with enveloping algebra $U(\g)$. The group $\aut (U(\g))$ of $\FF$-algebra automorphisms of $U(\g)$ is still for the most part unknown (except in particular instances, e.g.\ if $\g$ is abelian or if $\dim_\FF\g\leq 2$). For example, if 
$\g$ is the two-dimensional abelian Lie algebra, then $U(\g)$ is the polynomial algebra in two indeterminates $\FF[x_{1}, x_{2}]$, whose group of automorphisms is generated by the \emph{elementary} automorphisms of the form
\begin{equation*}
	x_{i}\mapsto \lambda x_{i}+f(x_{j}), \q x_{j}\mapsto x_{j} \q\q (i\neq j)
\end{equation*}
with $\lambda\in\FF^*$ and $f(x_{j})$ a polynomial in the variable $x_{j}$ (\cite{hJ42}, \cite{vdK53}). In contrast with this simple description, the conjecture that the polynomial algebra in three variables over 
$\FF$ has \emph{wild} automorphisms (i.e.\ automorphisms not in the subgroup generated by the elementary automorphisms) has been settled by Shestakov and Umirbaev (see~\cite{SU04}) assuming $\FF$ has characteristic $0$. Other examples are the enveloping algebra of $\mathfrak{sl}_{2}$ and the enveloping algebra of the Heisenberg Lie algebra $\h$, which are known to have wild automorphisms by results of Joseph~\cite{aJ76} and Alev~\cite{jA86}.

Concerning the Weyl algebra, its group of automorphisms was described by Dixmier in~\cite{jD68} in case $\chara(\FF)=0$ and later by Makar-Limanov in \cite{ML84} for arbitrary fields. We briefly describe this group next.

Let $\mathsf{SL_2(\FF)}$ denote the special linear group of $2 \times 2$ matrices over $\FF$ of
determinant $1$.   Each matrix  $\tt{S} =  \left( \begin{smallmatrix} \alpha & \gamma \\  \beta & \varepsilon \end{smallmatrix} \right)  \in \mathsf{SL_2(\FF)}$ determines an automorphism $\varphi_{\tt S}$  of $\wa(\FF)$ given by 

\begin{equation*}   x  \mapsto  \alpha x + \beta y,   \qquad    y \mapsto   \gamma x + \varepsilon y. \end{equation*}

The matrix $\tt T= \left( \begin{smallmatrix} 0 & 1 \\  -1 & 0 \end{smallmatrix} \right) \in \mathsf{SL_2(\FF)}$
corresponds to the automorphism $\tau = \varphi_{\tt T}$  of $\wa(\FF)$ given by  $x \mapsto -y$, \, $y \mapsto x$ 
and $\tau^{-1}$ corresponds to the automorphism with $x \mapsto y$, \, $y \mapsto -x$.   Note that $\tau^2 = -1$ and $\tau^4 = 1$.

For each $f \in \FF[x]$,  there is an automorphism $\phi_f$ with $\phi_f(x) = x$ and
$\phi_f(y) = y+f$. Let $\psi_f = \tau^{-1} \circ  \phi_{-f} \circ \tau$ and observe that
\begin{align*} \psi_f(x) = x + f(y)\q \text{and}\q
	\psi_f(y) =  y. \end{align*}

The following provide generating sets of automorphisms for $\aut(\wa(\FF))$ (compare 
\cite {ML84} and \cite{tS87},  and see also \cite{KA11} for part (iii)).

\begin{theorem}[{\cite[Thm.\ 8.19]{BLO15tams}}]\label{T:autgens}   Each of the following sets  gives a generating set for the automorphism
	group $\aut(\wa(\FF))$:
	\begin{itemize} 
		\item[{\rm (i)}] $\{\phi_f  \mid f \in \FF[x]\} \cup \{\psi_f \mid f \in \FF[x]\}$,
		\item[{\rm (ii)}]  $\{\varphi_{\tt S} \mid {\tt S} \in \mathsf{SL}_2(\FF) \}  \cup
		\{\phi_f  \mid f \in \FF[x]\}$,
		\item[{\rm (iii)}] $\{\tau, \phi_f  \mid f \in \FF[x]\}$,
		\item[{\rm (iv)}] $\{\tau, \psi_f  \mid f \in \FF[x]\}$.
	\end{itemize}
\end{theorem} 

Switching to derivations, it is well known that all derivations of $\wa(\FF)$ are inner (i.e., of the form $\ad_a$, for some $a\in\wa(\FF)$) in case $\chara(\FF) = 0$. In spite of a claim in \cite[Cor.]{pR73} that the same holds in positive characteristic, that statement is false, as shown in \cite{BLO15ja}.

Indeed, if $\chara(\FF)=p > 0$, the derivations $ (\ad_x)^p =\ad_{x^p}$ and $(\ad_y)^p =\ad_{y^p}$ are $0$ on the Weyl algebra  $\wa(\FF)$.   However, $\wa(\FF)$ has two special derivations
$E_x$ and $E_y$,  which are specified by
\begin{equation*} E_x(x) = y^{p-1},  \ \   E_x(y) = 0, \quad \hbox{and} \quad   E_y(x) = 0, \ \ E_y(y) = x^{p-1}.
\end{equation*}

\begin{theorem}[{\cite[Thm.\ 3.8]{BLO15ja}}]  Assume that $\chara(\FF) = p > 0$. 
	Then 
	\begin{itemize}
		\item[{\rm (a)}]  $\der(\wa(\FF)) =\mathsf{Z}(\wa(\FF)) E_x \oplus \mathsf{Z}(\wa(\FF)) E_y \oplus  \inder(\wa(\FF))$.
		\item[{\rm (b)}]  $\hoch(\wa(\FF)) =  \der(\wa(\FF))/\inder(\wa(\FF)) \cong \der(\FF[t_1, t_2])$ as Lie algebras, where $t_1 = x^p$, $t_2 = y^p$.  
	\end{itemize}  
\end{theorem} 

\subsection{The Dixmier and Jacobian conjectures}\label{SS:intro:wa:djc}

In \cite[Problem 1]{jD68}, Dixmier asked if every algebra endomorphism of $\wa[n](\FF)$ must be an automorphism when $\chara(\FF) = 0$. This is known as the Dixmier conjecture $DC_n$ and it is still open, even for $n=1$.

The Jacobian conjecture $JC_n$ is a statement in algebraic geometry that any polynomial endomorphism $\phi$ of the affine $n$-space $\FF^n=\mathrm{Spec} (\FF[x_1,\ldots, x_n])$ with Jacobian $1$, i.e.\ such that
\begin{align*}
	\det\seq{\frac{\partial \phi^*(x_i)}{\partial x_j}}_{1\leq i,j\leq n}=1,
\end{align*}
is an automorphism. The Jacobian conjecture holds trivially for $n=1$ but is otherwise open, even for $n=2$. 

The Dixmier conjecture $DC_n$ implies the Jacobian conjecture $JC_n$ and in \cite{BKK07} Belov-Kanel and Kontsevich show that $JC_{2n}$ implies $DC_n$, so in particular these two conjectures are stably equivalent (see also \cite{yT05}).

The Jacobian (and hence also the Dixmier) conjecture has many other equivalent forms, one of which concerns locally nilpotent derivations on polynomial rings and on Poisson algebras (see Subsection~\ref{SS:intro:wa:deform}). In fact, there seems to be a close connection in general between locally nilpotent derivations of (possibly noncommutative or nonassociative) algebras and automorphism groups, which has motivated our note in \cite{KLM21} connecting the theorems of Rentschler \cite{rR68} and Dixmier \cite{jD68} on locally nilpotent derivations and automorphisms of the polynomial ring $\FF[x_1, x_2]$ and of the Weyl algebra $\wa(\FF)$.

\subsection{The Weyl algebra: combinatorics}\label{SS:intro:wa:comb}

Surprising as it may be, there has been a great deal of combinatorics associated with the operators of multiplication by $t$ and $\frac{d}{dt}$ on $\FF[t]$, which can be traced back to Euler and Cayley (see for example \cite{hfS23} and \cite{aC57}). More recent variations can be found in \cite{lC73, BR87,mM11} and a comprehensive account of the subject is \cite{MS16}.

From a different perspective, let $\mathcal Y$ denote Young's lattice of partitions, i.e.\ the set of all integer partitions ordered by inclusion, and denote by $\FF\mathcal Y$ the vector space of formal linear combinations of elements of $\mathcal Y$. Then $\FF\mathcal Y$ affords a representation of the Weyl algebra $\wa(\FF)$, where $x$ acts by sending a partition $\lambda$ to the sum of all the partitions which can be obtained from $\lambda$ by adding $1$ to a part of $\lambda$ (without changing the property that the parts are weakly decreasing) or creating a new part $1$ at the end, and $y$ acts similarly but by subtracting $1$. Then, as operators on $\FF\mathcal Y$, we have $yx-xy=1$.

\begin{figure}[htbp]
	\begin{center}
		\scalebox{.7}{
			\xymatrix{{\yng(5)}\ar@{-}[dr]&&{\yng(4,1)} \ar@{-}[dl] \ar@{-}[dr] &{\yng(3,2)} \ar@{-}[d] \ar@{-}[dr]& {\yng(3,1,1)} \ar@{-}[dl] \ar@{-}[dr]  &{\yng(2,2,1)}   & {\yng(2,1,1,1)} \ar@{-}[dl] \ar@{-}[dr] &&{\yng(1,1,1,1,1)}\ar@{-}[dl]\\ 
				&{\yng(4)} \ar@{-}[dr] && {\yng(3,1)} \ar@{-}[dl] \ar@{-}[dr] & {\yng(2,2)} \ar@{-}[ur] \ar@{-}[d] & {\yng(2,1,1)} \ar@{-}[u] \ar@{-}[dl] \ar@{-}[dr] && {\yng(1,1,1,1)} \ar@{-}[dl]  \\ 
				&& {\yng(3)} \ar@{-}[dr] && {\yng(2,1)} \ar@{-}[dl] \ar@{-}[dr] && {\yng(1,1,1)} \ar@{-}[dl] \\ 
				&&& {\yng(2)} \ar@{-}[dr]  && {\yng(1,1)} \ar@{-}[dl] \\
				&&&& {\yng(1)} \ar@{-}[d]\\
				&&&& \emptyset} 
		}
		\caption{Young's lattice}
		\label{Ylattice}
	\end{center}
\end{figure}
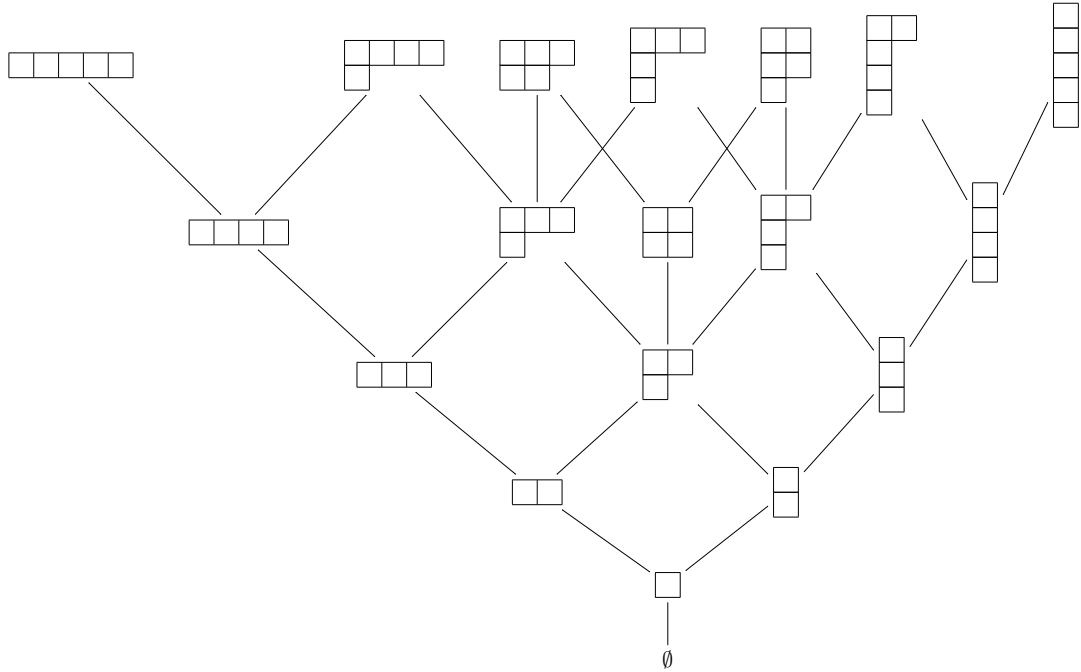

\begin{example}
	Representing partitions by their Young diagram, we have:
	\begin{align*}
		yx.\yng(2,1)&=y.\seq{\yng(3,1)+\yng(2,2)+\yng(2,1,1)}\\
		&=3\,\yng(2,1)+\yng(3)+\yng(1,1,1)\\[5pt]
		xy.\yng(2,1)&=x.\seq{\yng(1,1)+\yng(2)}\\
		&=2\,\yng(2,1)+\yng(3)+\yng(1,1,1).
	\end{align*}
	So
	\begin{align*}
		(yx-xy).\yng(2,1)=\yng(2,1).
	\end{align*}
\end{example}

In \cite{rS88}, Stanley showed that many of the properties of Young's lattice can be deduced from the commutation relation $yx-xy=1$ satisfied by the operators giving the action of the Weyl algebra on $\FF\mathcal Y$. This motivates the notion of  $r$-differential posets which are introduced and studied in \cite{rS88} and extended in \cite{rS90}, a theory which has interesting applications to the Sperner property on posets. A similar theory was independently developed by Fomin in \cite{sF86, sF94}.

Motivated by these ideas, Benkart and Roby introduced in \cite{BR98} a class of algebras named \textit{down-up algebras} which have been extensively studied from numerous points of view (see Definition~\ref{D:def:dua} and Subsection~\ref{SS:isoautder:gdua} for more details.). Our own work intersects these ideas in two directions:
\begin{enumerate}[label=\textup{(\roman*)}]
	\item We have introduced and studied a vast generalization of down-up algebras in \cite{LR22ca, LR22jaa};
	\item Computations stemming from work in representation theory initiated in \cite{BLO13} led us generalize a lot of the combinatorics associated to the Weyl algebra in \cite{BLR20}.
\end{enumerate}

\subsection{The Weyl algebra: ring theoretical properties}\label{SS:intro:wa:ring}

The Weyl algebra $\wa(\FF)$ has many interesting ring theoretical properties, some of which have already been mentioned. Firstly, it is a Noetherian domain of Gelfand-Kirillov dimension $2$ and it has a PBW-type basis $\pb{x^iy^j\mid i, j\geq 0}$. Both of these properties follow immediately form viewing it as an Ore extension. Since Ore extensions (also referred to as skew polynomial rings) underly most of our work, we briefly give a definition, also to establish notation.

A left Ore extension $R[\theta;\sigma, \delta]$  is built from a unital associative (not necessarily commutative) ring $R$, a ring endomorphism $\sigma$ of $R$ and a left $\sigma$-derivation of $R$, where by a  left $\sigma$-derivation $\delta$,  we mean that  $\delta$ is additive and $\delta (rs) = \delta(r)s + \sigma (r) \delta (s)$ holds for all $r,s \in R$. Then $R[\theta, \sigma, \delta]$ is the unital associative ring generated by $\theta$ over $R$ subject to the relation
\begin{align*}
	\theta r= \sigma(r)\theta + \delta(r) \qquad \hbox{for all} \ r \in R.
\end{align*}
Then it can be shown that $R[\theta;\sigma, \delta]$ is a free left $R$-module with free basis $\pb{\theta^i}_{i\geq 0}$, so that any element in $R[\theta;\sigma, \delta]$ can be uniquely expressed as a (noncommutative) polynomial in $\theta$ with coefficients on the left. Right Ore extensions $[\theta;\sigma, \partial]R$, are defined similarly for a right $\sigma$-derivation $\partial$ of $R$, in which case $[\theta;\sigma, \partial]R$ is a free right $R$-module with free basis $\pb{\theta^i}_{i\geq 0}$. If $\sigma$ is an automorphism of $R$ then both concepts are equivalent. Unless otherwise noted, we will consider left $\sigma$-derivations and left Ore extensions. In case $\sigma=1_R$, then it is customary to write simply $R[\theta; \delta]$ and similarly to write $R[\theta;\sigma]$ in case $\delta=0$.
The books \cite{GW89} and \cite{MR00} are excellent references for Ore extensions and Noetherian rings in general.

\begin{example}
	Let $\delta=\frac{d}{dx}$ be the ordinary derivation operator on $\FF[x]$. Then the Ore extension $\FF[x][y;\delta]$ is precisely the Weyl algebra $\wa(\FF)$. 
\end{example}

In case $\chara(\FF)=0$, then important additional properties of $\wa(\FF)$ are that it is simple, has trivial center equal to $\FF$, is rigid and has global dimension $1$. Given the discussion in Subsection~\ref{SS:intro:wa:djc} and the connection to the theory of $\mathcal D$-modules, representation theory and mathematical Physics, it is not surprising that there is a multitude of results and questions in ring theory concerning the Weyl algebras over fields of characteristic $0$. For the sake of brevity, we only refer here to the properties which most relate to our work. In case $\chara(\FF)=p>0$, then $\wa(\FF)$ is Azumaya over its center, has global dimension $2$ and is no longer rigid.

Another construction which has now become standard in ring theory and form which the Weyl algebra can be obtained is that of a generalized Weyl algebra (GWA, for short). Generalized Weyl algebras were introduced by Bavula in \cite{vB91} and are defined as follows. Given a ring $R$, an automorphism $\sigma$ of $R$ and a central element $a\in R$, the generalized Weyl algebra $R(\sigma, a)$ is the ring extension of $R$ generated by $x$ and $y$, subject to the relations:
\begin{equation}\label{E:noeth:gwa}
	yb=\sigma (b)y, \quad bx=x\sigma (b), \quad xy=a, \quad yx=\sigma (a)  \quad\quad \text{for all $b\in R$.}
\end{equation}
\begin{example}
	Let $R=\FF[t]$, $\sigma\in\aut(\FF[t])$ with $\sigma(t)=t+1$ and $a=t$. Then the generalized Weyl algebra $\FF[t](\sigma, t)$ is the Weyl algebra $\wa(\FF)$. 
\end{example}

\subsection{Deformation theory and Hochschild (co)homology}\label{SS:intro:wa:deform}

The idea of algebraic deformation parallels the theory of deformations of complex analytic structures, initiated in \cite{KS58}. 

A formal deformation of an associative algebra $A$ is a $\FF\formal{\hbar}$-algebra structure $A_\mu$ on $A\formal{\hbar}$:
\begin{align*}\label{E:intro:wa:ring:star}
	a\star b=ab+\mu_1(a, b)\hbar +\mu_2(a, b)\hbar^2+\mu_3(a, b)\hbar^3 + \cdots
\end{align*}
for all $a, b\in A$, where $ab$ is the product in $A$. Thus, we retrieve $A$ on setting $\hbar=0$. 

The associativity of $\star$ is controlled by the Hochschild cohomology of $A$, a (co)homology theory for associative algebras due to Gerstenhaber \cite{mG64}. We denote the $n$-th Hochschild homology and cohomology groups of $A$ by $\hoch_n(A)$ and $\hoch^n(A)$, respectively, and set
\begin{align*}
	\hoch_\bullet(A)=\bigoplus_{n\geq 0}\hoch_n(A),\q\q \hoch^\bullet(A)=\bigoplus_{n\geq 0}\hoch^n(A).
\end{align*}
Then $\hoch^0(A)=\mathsf{Z}(A)$ is the center of $A$ and $\hoch^1(A)=\der(A)/\inder(A)$ is the Lie algebra of outer derivations of $A$. The Hochschild cohomology $\hoch^{\bullet}(A)$ can be made into a Lie module for the Lie algebra $\hoch^1(A)$ under the so-called Gerstenhaber bracket. We remark that the Hochschild (co)homology is an invariant under derived equivalence.

Under this theory, and assuming that $\chara(\FF)=0$, we can see the Weyl algebra $\wa(\FF)$ as a deformation of the commutative polynomial algebra $\FF[x,y]$, under the Weyl--Groenewold product. Indeed, define for $a, b\in\FF[x,y]$,
\begin{align}
	a\star b=\sum_{n\geq 0} \frac{1}{n!}\frac{\partial^n a}{\partial y^n}\,\frac{\partial^n b}{\partial x^n}\,\hbar^n.
\end{align}
Then $\star$ defines an associative product on the vector space $\FF[x,y]\!\!\formal{\hbar}$ with
\begin{eqnarray*}
	x\star x =x^2, && y\star y=y^2,\\
	y\star x=yx+\hbar, && x\star y=xy.
\end{eqnarray*}
So 
\begin{align*}
	y\star x-x\star y=\hbar.
\end{align*}
Setting $\hbar=1$ we retrieve the Weyl algebra $\wa(\FF)$ as a deformation of the commutative polynomial algebra $\FF[x, y]$.

The above example is also related to a process called semiclassical limit. Let $A$ be a torsionfree $\mathcal{R}$-algebra and choose $0\neq h\in\mathcal{R}$ which is not a unit in $A$. Since $h$ is a central non-unit, $hA$ is a proper ideal of $A$. We will use any of the notations $\barr a$, $a+hA$ or the more suggestive $\eval{a}{h=0}$ to denote the image of $a\in A$ under the canonical map onto $\barr A = A/hA$. Note that the $\mathcal{R}$-algebra $\barr A$ is naturally an $\barr{\mathcal{R}}=\mathcal{R}/h\mathcal{R}$-algebra as well.

Assume that $\barr A$ as above is commutative. Then we can define a Poisson bracket on $\barr A$ by setting
\begin{equation}\label{E:def:semiclassical:poisson}
	\pb{\barr a, \barr b}=\barr{h^{-1}[a, b]}=\eval{(h^{-1}[a, b])}{h=0}, \quad \forall a, b\in A,
\end{equation}
where $h^{-1}[a, b]$ just denotes the unique element $\gamma(a, b)\in A$ such that $[a, b]=h\gamma(a, b)\in hA$ (the existence of such an element follows from the commutativity of $\barr{A}$ and the uniqueness from the fact that $h\neq 0$ is not a zero divisor in $A$). Indeed, it is straightforward to check that \eqref{E:def:semiclassical:poisson} is independent of the choice of representatives $a, b\in A$ and defines an $\barr{\mathcal{R}}$-bilinear Poisson bracket on $\barr{A}$. Endowed with this bracket, the Poisson algebra $\barr{A}$ is called the semiclassical limit of $A$ and, in turn, $A$ is called a quantization of $\barr{A}$. We refer to \cite[III.5.4]{BG02} and \cite[Sec.\ 1.1.3]{fD11} for more details and examples.

The algebra $A=\FF[x,y]\!\!\formal{\hbar}$ endowed with the product $\star$ defined in \eqref{E:intro:wa:ring:star}, along with $\mathcal R=\FF\formal{\hbar}$ and $h=\hbar$ satisfy the settings above with $\barr A=\FF[x,y]$. The Poisson algebra structure defined on $\barr A$ is the so-called Poisson-Weyl algebra with $\pb{y,x}=1$. It is not difficult to show that the Jacobian conjecture $JC_2$ (see Subsection~\ref{SS:intro:wa:djc}) is equivalent to the statement that every Poisson derivation $D$ of the Poisson-Weyl algebra such that $D(s)=1$ for some $s\in\barr A$ is locally nilpotent.

\section{Representation Theory}\label{S:rep}

\begin{flushright}{\it
		``Very roughly speaking, representation theory studies symmetry in linear spaces. It is a beautiful mathematical subject which has many applications, ranging from number theory and combinatorics to geometry, probability theory, quantum mechanics and quantum field theory.''} \\
	From \cite[Introduction]{pE11} by P.\ Etingof \textit{et al.}
\end{flushright}

\bigskip

Representation theory could be defined as a way to understand abstract objects through symmetry, i.e., through the ways in which they act as (general or structure-preserving) transformations on various spaces. In our work, we have investigated representation theory of several classes of algebras and we briefly recount this in this section. 

\subsection{Quantized enveloping algebras}\label{SS:rep:uq+}

Assume that $\chara(\FF)=0$ and fix a parameter $q\in\FF^*$ which is not a root of unity. Let $\g$ be a finite-dimensional complex semisimple Lie algebra. The quantized enveloping algebra associated with $\g$ is a one-parameter Hopf algebra deformation $U_q(\g)$ of the universal enveloping algebra $U(\g)$ of $\g$. Quantized enveloping algebras were defined independently by Drinfeld \cite{vD85} and Jimbo \cite{mJ85} in connection with the quantum Yang-Baxter equations and have since become objects of interest due to their relation with numerous areas of Mathematics and Physics.

Fix a triangular decomposition $\g=\g^+\oplus\h\oplus\g^-$ of $\g$, where $\h$ is a Cartan subalgebra and $\g^+$, $\g^-$ form a pair of maximal nilpotent subalgebras of $\g$. Accordingly, the algebra $U_q(\g)$ admits a triangular decomposition,
\[
U_q(\g)\simeq U_q(\g^+)\otimes U_{q}(\h)\otimes U_q(\g^-).
\] 
It seems natural to think of $U_q(\g^+)$ as a deformation of the enveloping algebra $U(\g^{+})$ of the nilpotent Lie algebra $\g^{+}$. 

Let $\g=\fsl_n$ be the complex semisimple Lie algebra of traceless $n\times n$ matrices and consider its maximal nilpotent subalgebra $\g^{+}=\fsl_n^{+}$ consisting of the strictly upper triangular matrices in $\g$. We use the standard notation $[k]=\frac{q^{k}-q^{-k}}{q-q^{-1}}$ for the (symmetric) $q$-version of the integer $k\in\ZZ$. 

\begin{definition}\label{D:uqsl+}
	The quantized enveloping algebra $U_q(\fsl_n^{+})$ is the associative unital $\FF$-algebra given by the Chevalley generators $e_1$, \ldots, $e_{n-1}$, subject to the so-called quantum Serre relations
	\begin{align}\label{E:bsu:def:qsr}
		e_{i}e_{j} - e_{j}e_{i} =&0 \q\q \text{if}\ \left| i-j \right| \neq 1,\\  \label{E:int:sl:rd:rel2}
		e_{i}^{2}e_{j} -(q + q^{-1})e_{i}e_{j}e_{i} + e_{j}e_{i}^{2}  =&0 \q\q \text{if}\  \left| i-j \right| = 1.
	\end{align} 
\end{definition}

\begin{example}
	In case $n=3$ the algebra $\fsl_3^{+}$ is just the Heisenberg Lie algebra and $U_q(\fsl_3^{+})$ can be presented by generators (different from the Chevalley generators) $X$, $Y$, $Z$, satisfying the relations:
	\begin{equation*}
		ZX=q^{-1}XZ, \q\q ZY=qYZ, \q\q XY-q^{-1}YX=Z.
	\end{equation*}
	This is known as the quantum Heisenberg algebra.
\end{example}

Lie's theorem implies that if $\FF$ is algebraically closed, then all finite-dimensional simple $U(\g^+)$-modules are $1$-dimensional. This result still holds for $U_q(\g^+)$ but the representation theories of $U(\g^+)$ and $U_q(\g^+)$ diverge when considering infinite-dimensional simple modules. 

We recall at this point that a (left) primitive ideal of a ring $R$ is the annihilator of a simple (left) $R$-module. Thus, simple $R$-modules with annihilator $I$ are in bijection with faithful simple $R/I$-modules. In general, any primitive ideal is prime and any maximal ideal is primitive. If $R$ is a commutative ring or a finite-dimensional algebra, then there is no distinction between maximal and primitive ideals but for noncommutative rings these concepts are different and there are even rings for which $(0)$ is a left but not right primitive ideal. 

Dixmier's claim that in general the representations of non-abelian Lie algebras defy classification has led researchers to focus first on classifying the primitive ideals of enveloping algebras, their quantizations and more generally of infinite-dimensional noncommutative algebras. The set of (left) primitive ideals of a ring $R$ is denoted by $\mathrm{Prim}(R)$ and it is a topological subspace of the space $\mathrm{Spec}(R)$ of prime ideals of $R$, endowed with the so-called Jacobson topology. Hence, the closed subsets of $\mathrm{Spec}(R)$ are those of the form 
$\{ P\in \mathrm{Spec}(R)\mid P\supseteq I \}$ for $I$ an ideal of $R$.

By \cite[Thm.\ 4.7.9]{jD96}, the primitive ideals of $U(\g^+)$ are maximal and the corresponding factor algebras are isomorphic to Weyl algebras. This no longer holds for $U_q(\g^+)$ as there are primitive ideals of this algebra which are not maximal; and even for the maximal ideals of $U_q(\g^+)$, it is no longer the case that the corresponding factor algebra (which will be simple) is a Weyl algebra. However, generalizing classical results of {So\u{\i}bel'man} for unitary representations of a maximal compact subgroup of $\mathrm{SL}(n)$, Hodges, Levasseur, Joseph and others have established a connection between primitive ideals of quantized enveloping algebras and symplectic leaves of Poisson structures on the corresponding algebraic groups.

\subsubsection{Primitive ideals of \texorpdfstring{$U_q(\fsl_n^+)$}{Uq(sln+)}}\label{SS:rep:uq+:prim}

By work of Alev and Dumas~\cite{AD94}, and Caldero~\cite{pC94a, pC95}, there exist homogeneous elements $\Delta_1, \ldots, \Delta_{n-1}$ in $U_q(\fsl_n^+)$ which commute with the Chevalley generators $e_1$, \ldots, $e_{n-1}$ up to a power of $q$ and which generate a polynomial algebra $\FF[\Delta_1, \ldots, \Delta_{n-1}]$. The center of $U_q(\fsl_n^+)$ is the polynomial algebra in the variables $z_i=\Delta_i\Delta_{n-i}$, for $i=1, \ldots, \lfloor (n-1)/2\rfloor$ and, in case $n=2k$, $z_k=\Delta_k$. So $\mathsf{Z}(U_q(\fsl_n^+))=\FF[z_1, \ldots, z_\ell]$, where $\ell=\lfloor n/2\rfloor$.

Given $\barr \alpha=(\alpha_1, \ldots, \alpha_\ell)\in\FF^{\ell}$, let
\begin{align*}
	I^{\mathbf{\barr\alpha}}=\sum_{i=1}^\ell U_q(\fsl_n^+)(z_i-\alpha_i),
\end{align*}
the ideal of $U_q(\fsl_n^+)$ generated by the central elements $z_i-\alpha_i$, for $i=1, \ldots, \ell$.

\begin{theorem}[{\cite[Thm.\ 3.5]{sL06}}]\label{T:rep:uq+:prim:mt}
	Assume $\barr \alpha\in(\FF^{*})^{\ell}$. Then $I^{\mathbf{\barr\alpha}}$ is a maximal ideal of $U_q(\fsl_n^+)$. It is also minimal among primitive ideals of $U_q(\fsl_n^+)$.
\end{theorem}

\begin{corollary}[{\cite[Cor.\ 3.6]{sL06}}]\label{C:rep:uq+:prim:gkd}
	Let $\barr \alpha\in(\FF^{*})^{\ell}$. The factor algebra $U_q(\fsl_n^+)/I^{\mathbf{\barr\alpha}}$ is a simple Noetherian domain with center $\FF$ and GK-dimension $\binom{n}{2}-\ell$. In particular, the GK-dimension of $U_q(\fsl_n^+)/I^{\mathbf{\barr\alpha}}$ is always even but $U_q(\fsl_n^+)/I^{\mathbf{\barr\alpha}}$ is not isomorphic to a Weyl algebra $\wa[k](\FF)$ for any $k\geq 1$.
\end{corollary}

The relevance of Corollary~\ref{C:rep:uq+:prim:gkd} is that it constructs explicitly new simple Noetherian domains with trivial center and even GK-dimension, which are not isomorphic to any Weyl algebra. Yet, they are analogues of the Weyl algebras in that they control part of the representation theory of the quantum nilpotent algebras $U_q(\fsl_n^+)$. At the time this result was published, only a few small examples of primitive quotients of quantum nilpotent enveloping algebras had been constructed, namely, in case $n=3$, \cite{tH90} (in connection with oscillator representations of quantized enveloping algebras), \cite{KS93} and \cite{mpM94}, where the corresponding factor algebra was named the Weyl-Hayashi algebra, and \cite{sL07} for the quantum nilpotent enveloping algebra $U_q(\mathfrak{so}_5^+)$, of type $B_2$.

\subsubsection{Infinite-dimensional representations of \texorpdfstring{$U_q(\fsl_n^+)$}{Uq(sln+)}}\label{SS:rep:uq+:reps}

A famous result of Kostant~\cite{bK63} (see also \cite[8.2.4]{jD96}), known as ``separation of variables'' shows that $U(\g)$ is free as a module over its center $\mathsf{Z}(\g)$. More concretely, he showed that multiplication gives an isomorphism $U(\g)\simeq\mathcal{H}\ot \mathsf{Z}(\g)$, where $\mathcal{H}$ is the so-called space of harmonics. The quantum analogue of this result for $U_q(\g)$, where $q$ is assumed to be transcendental over the rationals, was proved in \cite{JL94}.

In \cite[Thm.\ 2]{sL05}, assuming that $q$ is not a root of unity, we proved an analogous result for $U_q(\fsl_n^+)$. In fact, we showed the more general statement that $U_q(\fsl_n^+)$ is free as a module for the polynomial algebra $N=\FF[\Delta_1, \ldots, \Delta_{n-1}]$ defined in the previous subsection. Our motivation was the study of infinite-dimensional representations of $U_q(\fsl_n^+)$. Indeed, these results have important consequences for the representation theory of $U_q(\fsl_n^+)$, one of which being the existence of simple modules with arbitrary central character and, more generally, with arbitrary $N$-character. 

Indeed, we use the latter results to construct modules by inducing from one-dimensional $N$-modules. Given an 
$N$-character $\chi\in \widehat{N}$ with corresponding simple module 
$V_{\chi}=\FF v_{\chi}$, the induced $U_q(\fsl_n^+)$-module $M_{\chi}=U_q(\fsl_n^+)\otimes_{N}V_{\chi}$ has a weight space decomposition with respect to $N$,
\[
M_{\chi}=\bigoplus_{\eta\in \widehat{N}} M_{\chi}^{(\eta)},
\]
where $M_{\chi}^{(\eta)}=\{ m\in M_{\chi} \mid x.m=\eta (x)m \ \  \mbox{for all $x\in N$} \}$, and it is easy to see that every subquotient of $M_{\chi}$ inherits this grading. Thus, by considering maximal submodules of 
$M_{\chi}$, we obtain simple $U_q(\fsl_n^+)$-modules with nice properties, including a grading by the character group $\widehat{N}$. Specializing to central characters and assuming that $\FF$ is algebraically closed, the fact that $U_q(\fsl_n^+)$ satisfies the Nullstellensatz over $\FF$ (see \cite[Thm.\ 9.4.21]{MR00}) implies that all irreducible $U_q(\fsl_n^+)$-modules (even the infinite-dimensional ones) can be obtained as quotients of $M_{\chi}$, for an appropriately chosen $\chi\in \widehat{N}$.

Going back to the specific example of $U_q(\fsl_4^+)$, in \cite[Sec.\ 4.3]{sL06} we describe by generators all $4!=24$ primitive ideals of $U_q(\fsl_4^+)$, their heights and the GK-dimension of the corresponding factor algebras. For each primitive ideal, we constructed an irreducible representation of $U_q(\fsl_4^+)$ with that annihilator, thus in a way illustrating Dixmier's plan of first determining the primitive ideals and then, for each primitive quotient, studying the corresponding irreducible representations.

\subsection{The Taft algebra and its Drinfeld double}\label{SS:rep:ta}

Assume that $\FF$ is arbitrary and let $q\in \FF$ be a primitive $n$-th root of unity with $n \geq 2$. 

\begin{definition}
	The $n$-th Taft algebra is the unital associative $\FF$-algebra $U_n(q)$ with generators $G,X$ and relations
	\begin{align*}
		G^n = 1,\q X^n = 0, \q GX = qXG.
	\end{align*}
\end{definition}
This is an $n^2$-dimensional algebra with basis $\{ G^i X^j \mid 0\leq i,j\leq n-1 \}$.  It has a coalgebra structure given by $\Delta(G) = G \otimes G$, $\Delta(X) = X \otimes G + 1 \otimes X$, $\epsilon(G) = 1$, and $\epsilon(X) = 0$.  It also has a Hopf algebra antipode given by $S(G) = G^{-1}$ and $S(X) = -XG^{-1}$. 

The Drinfeld double $D(U_n(q))$ can be given in terms of generators $a,b,c,d$, subject to the relations
\begin{align*}
	a^n = 0 = d^n,\q  b^n=1 = c^n,\q  ba = qab,   \q  dc = qcd,\\ db = qbd,\q bc = cb,
	\q ca = qac, \q da - qad = 1 - bc.
\end{align*}
The Hopf algebra $D(U_n(q))$ and its quotient in the small quantum groups is further discussed in \cite[Chapter 9]{cK95}.
We view $U_n(q)$ as the Hopf subalgebra of $D(U_n(q))$ generated by $a,b$ via the embedding $G \mapsto b$ and $X \mapsto a$.

Let $Q$ denote the right $D(U_n(q))$-module $D(U_n(q))/U_n(q)^+D(U_n(q))$, where $U_n(q)^+$ is the augmentation ideal of $U_n(q)$. The module $Q$ is a generalization of the permutation module of cosets for a group extension. Using calculations on the Green ring of $U_n(q)$, in \cite{HKL17} we decomposed the $U_n(q)$-modules $Q$ and $Q\ot Q$ into a direct sum of indecomposables.

\begin{theorem}[\cite{HKL17}]\label{T:decQQ}
	All indecomposable $U_n(q)$-modules occur in the Krull-Schmidt decomposition of $Q\otimes Q$ as a right $U_n(q)$-module.
\end{theorem}

We remark that the endomorphism ring $\End Q^{\otimes 2}$ of the $U_n(q)$-module $Q \otimes Q$ is therefore Morita equivalent to the Auslander algebra of $U_n(q)$ (see \cite{SY11}). Extending several notions of depth for subalgebra extensions, Kadison defines in \cite{lK14} the notion of depth of a nonzero object by means of similarity relations between consecutive truncated tensor powers (a Morita invariant in terms of Morita invariance of ring extensions, see \cite{lK14, HKL17} for more details). Then we have the following.

\begin{corollary}[\cite{HKL17}]
	For the Taft algebra $U_n(q)$ in its Drinfeld double, the depth of $Q$ is $d(Q_{U_n(q)}) = 2$ and the minimum even depth $$d_{ev}(U_n(q), D(U_n(q)) = 6.$$
\end{corollary}

\subsection{The quantum plane and the quantum Weyl Algebra}\label{SS:rep:qpqwa}

Fix $q\in\FF^*$ and consider the operators $\tau_q$ and $\partial_q$ defined on the polynomial algebra $\FF[t]$ by
\begin{equation}\label{E:rep:qprep}
	\tau_q (p)(t)=p(qt), \quad \mbox{and} \quad \partial_q (p)(t)=\frac{p(qt)-p(t)}{qt-t}, \quad \mbox{for $p\in\FF[t]$}.
\end{equation}
These operators satisfy the relation $\partial_q\tau_q=q \tau_q\partial_q$ and if $q$ is not a root of unity, this relation defines the unital associative algebra which they generate.

\begin{definition}
	Let $\FF$ be a field, and let $q\in\FF^*$ with $q\neq 1$.  The quantum plane 
	is the unital associative algebra 
	\begin{equation}\label{E:rep:qp:def}
		\FF_q[x, y]=\FF\ip{ x, y}/(yx-qxy)
	\end{equation}
	with generators $x$ and $y$ subject to the relation $yx=qxy$. 
\end{definition}

Thus, the assignment $x\mapsto \tau_q$, $y\mapsto \partial_q$ yields a (reducible) representation of $\FF_q[x, y]$. The operators $\tau_q$ and $\partial_q$ are central in the theory of linear $q$-difference equations and $\partial_q$ is also known as the Jackson derivative, as it appears in \cite{fJ10}. See e.g.\ \cite{yM88}, \cite[Chap.\ IV]{cK95} and references therein for further details.

If in \eqref{E:rep:qprep} we replace the operator $\tau_q$ with the operator $\nu (p)(t)=tp(t)$, then instead we obtain the relation $\partial_q\nu=q\nu\partial_q + 1$.

\begin{definition}
	Let $\FF$ be a field, and let $q\in\FF^*$ with $q\neq 1$.  The (first) quantum Weyl algebra (or $q$-Weyl algebra)
	is the unital associative algebra 
	\begin{equation}
		\wa^q(\FF)=\FF\ip{ x, y}/(yx-qxy-1)
	\end{equation}
	with generators $x$ and $y$ subject to the relation $yx=qxy+1$. 
\end{definition}
The quantum Weyl algebra and its higher degree analogues are amongst the simplest non-abelian examples of quantum multiplicative quiver varieties (see the recent paper~\cite{CGJ20}) and have been extensively studied.

The irreducible representations of the quantum plane $\FF_q[x, y]$ and of the quantum Weyl algebra $\wa^q(\FF)$ have been classified in~\cite{vB97} using results from~\cite{BvO97}. Following \cite{vB97} we say that a representation of $\FF_q[x, y]$ is a \emph{weight representation} if it is semisimple as a representation of the polynomial subalgebra $\FF[H]$ generated by the element $H=xy$. When $q$ is a root of unity all irreducible representations of $\FF_q[x, y]$ are finite-dimensional weight representations, and these are well understood. For example, if $\FF$ is algebraically closed and $q$ is a primitive $n$-th root of unity then the irreducible representations of $\FF_q[x, y]$ are either $1$ or $n$ dimensional. When $q$ is not a root of unity there are irreducible representations of $\FF_q[x, y]$ that are not weight representations, and in particular are not finite dimensional. These turn out to be the \emph{$\FF[H]$-torsionfree} irreducible representations of $\FF_q[x, y]$, as they remain irreducible (i.e. nonzero) upon localizing at the nonzero elements of $\FF[H]$.  The torsionfree representations of $\FF_q[x, y]$ are classified in terms of elements satisfying certain conditions in~\cite[Cor.\ 3.3]{vB97}, but no explicit construction of these representations is given.

In \cite{LL15}, as part of an undergraduate research project with the student João Lourenço (now a research associate at Imperial College), we gave an explicit construction of a 3-parameter family $\mathsf{V}^{m,n}_{f}$ of infinite-dimensional representations of $\FF_q[x, y]$ having the following properties:
\begin{itemize}
	\item $m$ and $n$ are positive integers, and $f:\mathbb{Z} \rightarrow \FF^{*}$ satisfies the condition $f(i+n)=qf(i)$, for all $i\in\ZZ$, which essentially encodes $n$ independent parameters from $\FF^*$;
	\item $\mathsf{V}^{m,n}_{f}$ is irreducible if and only if $\gcd(m, n)=1$;
	\item if $(m, n)\neq (m', n')$ then $\mathsf{V}^{m,n}_{f}$ and $\mathsf{V}^{m',n'}_{f'}$ are not isomorphic;
	\item $\mathsf{V}^{m,n}_{f}$ is a weight representation if and only if $m=n$;
	\item if $\FF$ is algebraically closed and $V$ is an irreducible weight representation of $\FF_q[x, y]$ that is infinite dimensional, then $V\simeq\mathsf{V}^{1,1}_{f}$ for some $f:\mathbb{Z} \rightarrow \FF^{*}$.
\end{itemize}
Thus, in some sense weight and non-weight representations of $\FF_q[x, y]$ are rejoined in the family $\mathsf{V}^{m,n}_{f}$.

The localization of $\FF_q[x, y]$ at the multiplicative set generated by $x$ contains a copy of the quantum Weyl algebra $\wa^q(\FF)$. This is used to regard the representations $\mathsf{V}^{m,n}_{f}$ as infinite-dimensional irreducible representations of $\wa^q(\FF)$. In contrast with the action of $\FF_q[x, y]$ on $\mathsf{V}^{m,n}_{f}$ when $m=n$, it turns out that $\mathsf{V}^{m,n}_{f}$ is never a weight representation of $\wa^q(\FF)$ in the sense of~\cite{vB97}.

\subsection{A Parametric Family of Subalgebras of the Weyl Algebra}\label{SS:rep:ah}

Over a series of papers with G.\ Benkart and M.\ Ondrus (\cite{BLO15tams,BLO13,BLO15ja}), we defined and studied a family of infinite-dimensional unital associative algebras $\A_h$ parametrized by a polynomial $h
= h(x) \in \FF[x]$,
where $\FF$ is an arbitrary field.

\begin{definition}\label{D:rep:ah}
	Let $\FF$ be a field, and let $h \in \FF[x]$.  The algebra  $\A_h$ is  the unital associative algebra over $\FF$ with generators $x$, $y$ and defining relation $yx = xy + h$ (equivalently,  $[y,x] = h$). 
\end{definition}

These algebras appear naturally in the study of Ore extensions, due to the observation that any Ore extension of the form $\FF[x][y,\sigma,\delta]$, where $\sigma$ is an automorphism, is either a quantum plane $\FF_q[x, y]$, a quantum Weyl algebra $\wa^q(\FF)$ or an algebra in the family $\A_h$, with $h = h(x) \in \FF[x]$. Quantum planes and quantum Weyl algebras have been extensively studied and are somewhat more manageable (see in particular Subsection~\ref{SS:rep:qpqwa}), whereas the algebras $\A_h$ have many more parameters of freedom and are more closely connected with the Weyl algebra, in case $h\neq 0$.

\begin{example}
	Assume that $\chara(\FF)=0$. Replacing the differentiation operator on $\FF[t]$ by integration one obtains the operators 
	\begin{align*}
		X.p(t)=\int_{0}^{t}p(z)\, dz \q\q Y.p(t)=t\, p(t)
	\end{align*}
	and it follows that
	\begin{align*}
		XY.p(t)=\int_{0}^{t}z\,p(z)\, dz =t\int_{0}^{t}p(z)\, dz - \int_{0}^{t}\int_{0}^{z}p(w)\, dw\, dz=(YX-X^{2}).p(t).
	\end{align*}
\end{example}
The relation $XY-YX=X^2$ on the free algebra $\FF\ip{X,Y}$ defines the so-called Jordan plane, which in the notation above is the algebra $\A_{x^2}$. The Jordan plane arises in noncommutative algebraic geometry (see, for example, \cite{SZ94} and \cite{AS95})  and
exhibits many interesting features such as being Artin-Schelter regular of dimension 2. 
In a series of articles \cite{eS05, eS07, eS07b},  Shirikov
has undertaken an extensive study of the derivations, prime ideals, and modules
of the algebra $\A_{x^2}$.  These investigations have been extended by Iyudu \cite{nI14} to include results on finite-dimensional modules and automorphisms of
$\A_{x^2}$ over algebraically closed fields of characteristic zero.  Cibils, Lauve, and Witherspoon \cite{CLW09}  have used quotients of the algebra $\A_{x^2}$ 
and cyclic subgroups of their automorphism groups to
construct new examples of finite-dimensional Hopf algebras in prime characteristic
which are Nichols algebras.

Deformations of the Jordan plane, which can be seen as $\A_{x^2+\Delta}$, where $\Delta\in\FF$, have also appeared in many contexts including related to non-commutative probability theory~\cite{BKS97} and one-dimensional asymmetric exclusion models~\cite{DEHP93}.

Other striking examples of algebras in the family $\A_h$ are the following:
\begin{itemize}
	\item $\A_{0}=\FF[x, y]$, the (commutative) polynomial algebra;
	\item $\A_{1}=\wa(\FF)$, the Weyl algebra;
	\item $\A_{x}=U(\mathfrak{L})$, the enveloping algebra of the two-dimensional non-abelian Lie algebra 
	$\mathfrak{L}=\FF x\oplus\FF y$, with $[y, x]=x$.
\end{itemize}
There are striking similarities in the behavior of the algebras $\A_h$ as $h$ ranges over the polynomials
in $\FF[x]$ and that motivated our work. For example, if $0\neq f, g\in\FF[x]$ and $f$ divides $g$, then  $\A_g$ embedds into $\A_f$. In particular, $\A_h$ can be seen as a subalgebra of the Weyl algebra $\A_1$ for all $0\neq h\in\FF[x]$.

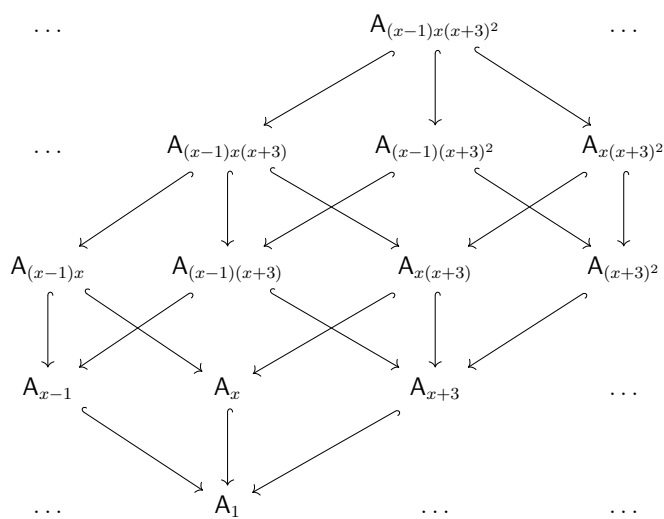
\begin{figure}[htbp]
	\begin{center}
		\scalebox{.8}{
			\begin{tikzcd}[ampersand replacement=\&]
				\ldots  \&  \& \A_{(x-1)x(x+3)^2} \arrow[hook]{dl} \arrow[hook]{d} \arrow[hook]{dr} \& \ldots \\[15pt]
				\ldots   \& \A_{(x-1)x(x+3)} \arrow[hook]{dl} \arrow[hook]{d} \arrow[hook]{dr} \& \A_{(x-1)(x+3)^2} \arrow[hook]{dl} \arrow[hook]{dr} \& \A_{x(x+3)^2} \arrow[hook]{dl} \arrow[hook]{d} \\[15pt]
				\A_{(x-1)x} \arrow[hook]{d} \arrow[hook]{dr}  \& \A_{(x-1)(x+3)} \arrow[hook]{dl}  \arrow[hook]{dr} \& \A_{x(x+3)} \arrow[hook]{dl} \arrow[hook]{d}  \& \A_{(x+3)^2} \arrow[hook]{dl}  \\[15pt]
				\A_{x-1}  \arrow[hook]{dr}  \& \A_{x} \arrow[hook]{d} \& \A_{x+3} \arrow[hook]{dl} \& \ldots \\[15pt]
				\ldots       \& \A_{1} \& \ldots \&  \ldots
			\end{tikzcd}
		}
		\caption{The lattice of inclusions in the family $\A_h$}
		\label{inclusions}
	\end{center}
\end{figure}

The representation theory of the family of algebras $\A_h$ over arbitrary fields was investigated in \cite{BLO13}, assuming that $h\neq 0$. The main results are descriptions of the irreducible modules, generalized weight $\A_h$-modules, and the primitive ideals of $\A_h$, which were shown to be either maximal ideals of finite co-dimension or the zero ideal (the latter in case $\chara(\FF)=0$). Results on indecomposable modules were also obtained in arbitrary characteristic.

If $\chara(\FF)>0$, all irreducible $\A_h$-modules are finite-dimensional and are factors of the induced modules $\A_h\otimes_{\FF[x]}\FF[x]/\mathfrak{m}$ for maximal ideals $\mathfrak m$ of $\FF[x]$. In case $\FF=\barr\FF$ then detailed descriptions were obtained.

In case $\chara(\FF)=0$, there are, in general, still finite-dimensional irreducible $\A_h$-modules but also infinite dimensional ones. Some are induced modules of the form $\A_h\otimes_{\FF[x]}\FF[x]/(f)$, where $f$ does not divide $h$, and others are $\FF[x]$-torsionfree irreducible modules. Using the embedding $\A_h\subseteq\A_1$, the latter were studied from the point of view of $\FF[x]$-torsionfree irreducible modules for the Weyl algebra $\A_1$.

\subsection{Quantum generalized Heisenberg algebras}\label{SS:rep:qgha}

In \cite{LR22ca} we introduced a new class of algebras, which we named \textit{quantum generalized Heisenberg algebras} (qGHA for short), as they can be seen simultaneously as deformations and as generalizations of the generalized Heisenberg algebras appearing in \cite{CRM01} and profusely studied thenceforth in the Physics literature (see e.g.~\cite{CRMR13}, \cite{BBH11}, \cite{BCG18} and the references therein). In the Mathematics literature, generalized Heisenberg algebras were studied mainly in \cite{LZ15}, \cite{LMZ15} and \cite{sL17}. For an overview of their relevance in mathematical Physics see the introductory section in \cite{LZ15}. 

\begin{definition}\label{D:qgha}
	Fix $q\in\FF$ and $f,g\in\FF[h]$.  The \textit{quantum generalized Heisenberg algebra} (qGHA, for short), denoted by $\hqfg$, is the $\FF$-algebra generated by $x$, $y$ and $h$, with defining relations:
	\begin{equation}\label{E:intro:def:qgha}
		hx=xf(h), \quad yh=f(h)y, \quad yx-qxy=g(h).
	\end{equation} 
\end{definition}

In \cite{LZ15}, working over the complex field $\CC$ and motivated by the Physics literature, the generalized Heisenberg algebra $\hh(f)$, parametrized by $f\in\CC[h]$, was introduced as the unital associative $\CC$-algebra with generators $x$, $y$ and $h$, with defining relations:
\begin{equation*}
	hx=xf(h), \quad yh=f(h)y, \quad yx-xy=f(h)-h.
\end{equation*}
It follows immediately that generalized Heisenberg algebras are precisely the qGHA (over $\CC$) with $q=1$ and $g=f(h)-h$, i.e.\ $\hh(f)=\hh_1(f, f-h)$. Working over an arbitrary field $\FF$, we can thus view qGHA as a generalization of the latter class of algebras, by deforming and generalizing the relation $yx-xy=f(h)-h$, turning it into a skew-commutation relation and allowing the skew-commutator to equal a generic polynomial, independent from $f$.

Besides providing a broader framework for the investigation of the underlying physical systems, our main motivation for introducing a generalization of this class comes from the observation in \cite{sL17} that the classes of generalized Heisenberg algebras and of (generalized) down-up algebras (see Definitions~\ref{D:def:dua} and~\ref{D:def:gdua} in Subsection~\ref{SS:isoautder:gdua}) intersect, although neither one contains the other. 

In spite of the name, the class of generalized Heisenberg algebras does not include the enveloping algebra of the Heisenberg Lie algebra nor its quantum deformations introduced in~\cite{KS93}, nor the enveloping algebra of $\mathfrak{sl}_2$. These and many other like algebras are now included in the class of qGHA and we show that they can all be studied uniformly with generalized Heisenberg algebras, highlighting their common properties.

The other interesting feature of our study comes from the fact that quantum generalized Heisenberg algebras are generically non-Noetherian and we have not found enough studies into the representation theory of non-Noetherian algebras which are nonetheless related to deformations of enveloping algebras of Lie algebras, as seems to be the case with qGHA. 

\subsubsection{Raising and lowering operators}\label{SS:intro:raiselower}

Consider the $3$-dimensional Lie algebra $\mathfrak{sl}_2$, with basis elements $x, y, h$ and Lie bracket given by $[x,h]=2x$, $[h,y]=2y$ and $[y,x]=h$. We can view its enveloping algebra as the qGHA $\hh_1(h-2, h)$. In the representation theory of $\mathfrak{sl}_2$, $x$ and $y$ are often represented as  raising and lowering operators on a finite or countable vector space. For example, in \cite{P82} and \cite{rS88} the existence of such operators on the vector space whose distinguished basis is a suitably defined poset is used to solve important combinatorial problems.

Take $V=\FF[t^{\pm 1}]$, the Laurent polynomial algebra, and suppose that $x$ and $y$ act on $V$ as raising and lowering operators, respectively, so that $h$ acts diagonally. We can assume that, relative to the basis $\left\{t^k\right\}_{k\in\ZZ}$, we have
\begin{equation}\label{E:intro:raislowLaurent}
	x t^k=t^{k+1},\quad  yt^k=\mu(k)t^{k-1}\quad  \mbox{and}\quad  ht^k=\lambda(k)t^{k},\qquad\mbox{for all $k\in\ZZ$,}
\end{equation}
where $\lambda, \mu:\ZZ\longrightarrow\FF$. Then, the $\mathfrak{sl}_2$ relations impose the conditions 
\begin{equation*}
	\mu(k+1)=\mu(k)+\lambda(k) \quad  \mbox{and}\quad \lambda(k+1)=\lambda(k)-2,
\end{equation*}
so $\lambda(k+1)$ is affine on $\lambda(k)$ and $\mu(k+1)$ is linear on $\mu(k)$ and $\lambda(k)$. 

Consider now the $3$-dimensional Heisenberg Lie algebra with basis elements $x, y, h$ and Lie brackets $[h,x]=[h,y]=0$ and $[y,x]=h$.
Its enveloping algebra can be seen as the qGHA $\hh_1(h, h)$. Then, the Heisenberg relations imposed on \eqref{E:intro:raislowLaurent} give
\begin{equation*}
	\mu(k+1)=\mu(k)+\lambda(k) \quad  \mbox{and}\quad \lambda(k+1)=\lambda(k),
\end{equation*}
so $\lambda$ is constant and $\mu(k+1)$ is affine on $\mu(k)$.

Another related example is given by the algebras similar to the enveloping algebra of $\mathfrak{sl}_2$ introduced by Smith in \cite{spS90}. These are precisely the qGHA of the form $\hh_1(h-1,g)$, for $g\in\FF[h]$. The corresponding conditions imposed on \eqref{E:intro:raislowLaurent} by the Smith algebra relations are
\begin{equation*}
	\mu(k+1)=\mu(k)+g(\lambda(k)) \quad  \mbox{and}\quad \lambda(k+1)=\lambda(k)-1,
\end{equation*}
so  $\lambda(k+1)$ is affine on $\lambda(k)$ but now $\mu(k+1)-\mu(k)$ is polynomial on $\lambda(k)$.

As a final example, if we take a generalized Heisenberg algebra, i.e.\ a qGHA of the form $\hh_1(f,f-h)$, then the corresponding conditions imposed on \eqref{E:intro:raislowLaurent} are
\begin{equation*}
	\mu(k+1)-\mu(k)=\lambda(k+1)-\lambda(k) \quad  \mbox{and}\quad \lambda(k+1)=f(\lambda(k)),
\end{equation*}
so $\lambda$ and $\mu$ differ by a constant and $\lambda(k+1)$ is polynomial in $\lambda(k)$.

With the more general relations allowed for by our definition of a qGHA, we can include all of the above cases and, more generally, in $\hqfg$ we have
\begin{equation*}
	\mu(k+1)=q\mu(k)+g(\lambda(k)) \quad  \mbox{and}\quad \lambda(k+1)=f(\lambda(k)),
\end{equation*}
so that $\lambda(k+1)$ is polynomial in $\lambda(k)$ and $\mu(k+1)$ is affine in $\mu(k)$ and polynomial in $\lambda(k)$. Representations of the qGHA $\hqfg$ will thus classify the creation and annihilation operators as in \eqref{E:intro:raislowLaurent}, under the latter assumptions.

\subsubsection{The finite-dimensional simple \texorpdfstring{$\hqfg$}{Hq(g,f)}-modules}\label{SSS:irreps}
In \cite{LR22ca} we completely classify the finite-dimensional simple $\hqfg$-modules for all polynomials $f, g\in\FF[h]$, assuming only that $q\neq 0$ and $\FF$ is algebraically closed. In particular, this study generalizes and unifies the classification of finite-dimensional simple modules over down-up algebras, generalized down-up algebras and generalized Heisenberg algebras, which has been carried out over the series of references \cite{BR98}, \cite{CM00}, \cite{CS04}, \cite{LZ15}, \cite{spS90}, \cite{lLB95} and \cite{sR02}, to name a few. 

The main result in \cite{LR22ca} is the following theorem. For the sake of brevity, we omit some of the definitions appearing below and refer the reader to \cite[Sec.\ 3]{LR22ca} for complete details. But to give an idea of the concepts involved we offer the following definitions.

\begin{itemize}
	\item $\displaystyle \Sf=\pb{\lambda:\ZZ\longrightarrow\FF\mid f(\lambda(i))=\lambda(i+1), \ \mbox{for all $i\in\ZZ$}}$.
	\item For $\lambda\in\Sf$, define $|\lambda|\geq 0$ so that  $\displaystyle |\lambda|\,\ZZ=\pb{k\in\ZZ\mid \lambda(i)=\lambda(i+k), \ \mbox{for all $i\in\ZZ$}}$.
	
	\item Given $\lambda\in\Sf$, set $\displaystyle \Tf=\pb{\mu:\ZZ\longrightarrow\FF\mid \mu(i+1)=q\mu(i)+g(\lambda(i)), \ \mbox{for all $i\in\ZZ$}}$.
	
	\item For $\lambda\in\Sf$ and $\mu\in\Tf$, define $|\mu|\geq 0$ so that $\displaystyle |\mu|\,\ZZ= \pb{k\in\ZZ\mid\mu(k|\lambda|)=\mu(0)}$.
	
	\item For $\lambda\in\Sf$ and $\mu\in\Tf$, define the $\hqfg$-module $\Af$ by
	defining the following action in $\Af=\FF[t^{\pm 1}]$
	\begin{equation*}
		ht^i=\lambda(i)t^i, \quad xt^i=t^{i+1}, \quad yt^i=\mu(i)t^{i-1}, \quad \mbox{for all $i\in\ZZ$.}
	\end{equation*}
	The $\hqfg$-module $\Bf$ is defined in a dual fashion.
\end{itemize}

\begin{theorem}\label{T:irreps:fdiml:all}
	Assume $\FF=\overline\FF$ and $q\neq 0$. Then any simple $n$-dimensional $\hqfg$-module is isomorphic to exactly one of the following simple modules:
	\begin{enumerate}[label=(\alph*)]
		\item $\Af/\FF[t^{\pm1}](t^{|\lambda||\mu|}-\gamma)$, for some $\lambda\in\Sf$, $\mu\in\Tf$ and $\gamma\in\FF^*$ such that $n=|\lambda||\mu|$.\label{T:irreps:fdiml:all:a}
		\item $\Bf/\FF[t^{\pm1}](t^{|\lambda||\mu|}-\gamma)$, for some $\lambda\in\Sf$, $\mu\in\Tf$ and $\gamma\in\FF^*$ such that $n=|\lambda||\mu|$ and $\mu(i)=0$ for some $0\leq i<|\lambda||\mu|$.\label{T:irreps:fdiml:all:b}
		\item $\Cf/\FF[t]t^n$, for some $\alpha\in\FF$ such that 
		$\nu_\alpha(i)\neq0$ for all $1\leq i\leq n-1$ and $\nu_\alpha(n)=0$.\label{T:irreps:fdiml:all:c}
	\end{enumerate}
\end{theorem}

\section{Automorphisms, isomorphisms and derivations}\label{S:isoautder}

It is quite common for different sets of generators and relations to yield the same intrinsic structure and this leads to classifying algebras by isomorphism. 

\begin{example}
	In \cite{DP09} the authors defined the \textit{$q$-meromorphic Weyl algebra} as the unital associative algebra with generators $x, y$ and defining relation $yx=qxy+x^2$, which appears to be a $q$-version of the Jordan plane $\A_{x^2}$ (see Definition~\ref{D:rep:ah}). But using the generator $Y=y+(q-1)^{-1}x$ instead of $y$ we find that $Yx=qxY$ so in fact the $q$-meromorphic Weyl algebra is the quantum plane $\FF_q[x, y]$ (see~\eqref{E:rep:qp:def}) in disguise.
\end{example}

The above example shows the importance of classifying a given family of algebras by isomorphism since seemingly distinct algebras can be essentially the same.

Automorphisms reflect the inner symmetries of an algebra and are thus an extremely useful tool for understanding it intrinsically. They also measure in a certain sense the rigidity of the algebra and to which extent a given set of generators and relations are unique. The interest in automorphism groups goes back to Cayley and is pivotal in Hilbert and Noether's work on invariant theory (consider for instance the theorem of Skolem and Noether on automorphisms of simple rings). Moreover, as illustrated in Subsection~\ref{SS:intro:wa:djc}, there are several important open problems related to automorphisms, even for (infinite-dimensional) algebras with as little as two generators.

There are several connections between automorphisms and derivations of algebras. For a finite-dimensional algebra over the complex field, its derivations form the Lie algebra of its automorphism group but such a close connection does not hold in general for infinite-dimensional algebras or algebras over fields of positive characteristic. 

However, over a field of characteristic $0$, from a locally nilpotent derivation $\partial$ of an algebra $A$ one can form the map $e^\partial=\sum_{k\geq 0}\frac{\partial^k}{k!}$, which is a well-defined automorphism of $A$. Moreover, gradings of $A$ by the base field $\FF$ are in correspondence with diagonalizable (hence locally finite) derivations of $A$. A final point of interest is that the first Hochschild cohomology $\hoch^1(A)$ consists precisely of the (outer) derivations of $A$.

\subsection{The case of \texorpdfstring{$U_q(\mathfrak{sl}^{+}_{4})$}{Uq(sl4+)} and the Andruskiewitsch-Dumas conjecture}\label{SS:isoautder:ADc}

In their paper~\cite{AD08}, Andruskiewitsch and Dumas conjectured that, given a finite\--di\-men\-sional complex simple Lie algebra $\mathfrak{g}$ with triangular decomposition 
$\mathfrak{g}=\mathfrak{g}^{-}\oplus\mathfrak{h}\oplus\mathfrak{g}^{+}$, then 
$\aut (U_{q}(\mathfrak{g}^{+}))$, the group of algebra automorphisms of the quantized enveloping algebra of the nilpotent Lie algebra $\mathfrak{g}^{+}$, is isomorphic to the semi-direct product of the torus $(\FF^{*})^n$ ($n$ being the rank of $\mathfrak{g}$) with the group of order $1$, $2$ or $3$ consisting of the diagram automorphisms of $\mathfrak{g}^{+}$, see~\cite[Prob.\ 1]{AD08}. This conjecture was first shown to hold for $\mathfrak{g}^{+}=\mathfrak{sl}^{+}_{3}$ in \cite{pC95, AD96} and for $\mathfrak{g}^{+}=\mathfrak{so}^{+}_{5}$ in \cite{sL07}.

In \cite{LL07} we computed the Lie algebra of derivations of $U_q(\mathfrak{sl}^{+}_{4})$ and showed that the first Hochschild cohomology group $\hoch^1(U_q(\mathfrak{sl}^{+}_{4}))$ is a free module of rank $3$ over the center of $U_{q}(\mathfrak{sl}^{+}_{4})$. To do this, we first apply the deleting derivations algorithm of Cauchon~\cite{gC03} so that, after suitably localizing, we can embed $U_{q}(\mathfrak{sl}^{+}_{4})$ in a quantum torus $T$. Extending a derivation $D$ of $U_{q}(\mathfrak{sl}^{+}_{4})$ to $T$ we obtain, by a result of Osborn and Passman~\cite{OP95}, a decomposition
\begin{equation*}
	D=\mathrm{ad}_{x}+\theta
\end{equation*}
with $x\in T$ and $\theta$ a central derivation of $T$. Using a sort of restoring derivations algorithm, we finish by deducing that $x\in U_{q}(\mathfrak{sl}^{+}_{4})$ and that $\theta$ sends each Chevalley generator of $U_{q}(\mathfrak{sl}^{+}_{4})$ to a multiple of itself by a central element of $U_{q}(\mathfrak{sl}^{+}_{4})$.

\begin{theorem}[{\cite[Thm.\ 3.11]{LL07}}]\label{T:aut:aut:mt:ders}
	$\hoch^1(U_q(\mathfrak{sl}^{+}_{4}))$ is free of rank $3$ over the center of $U_q(\mathfrak{sl}^{+}_{4})$. 
\end{theorem}

As an application, using the methods of~\cite{AC92} and~\cite{LL07jncg}, we confirmed the Andruskiewitsch-Dumas conjecture for $U_q(\mathfrak{sl}^{+}_{4})$.

\begin{theorem}[{\cite[Thm.\ 2.7]{LL07}}]\label{T:aut:aut:mt}
	$\aut (U_q(\mathfrak{sl}^{+}_{4}))$ is isomorphic to the semi-direct product of the $3$-torus $(\FF^*)^3$ and the group of order $2$ generated by the diagram  automorphism $\eta$ of $U_q(\mathfrak{sl}^{+}_{4})$.
\end{theorem}

In \cite{mY14}, the Andruskiewitsch-Dumas conjecture was proved in general.

\subsection{The case of \texorpdfstring{$\A_h$}{Ah}}\label{SS:isoautder:Ah}

In \cite{BLO15tams} we obtained the following classification results (recall Definition~\ref{D:rep:ah}).

\begin{theorem}[Isomorphism problem for $\A_h$ {\cite[Thm.\ 8.2]{BLO15tams}}]
	\begin{align*}
		\A_{h}\cong \A_{g}\iff  g(x)= \nu h(\alpha x+\beta),
	\end{align*}
	for some $\alpha, \beta, \nu\in\FF, \ \alpha\nu\neq 0$. 
\end{theorem}

Recall the definition of a generalized Weyl algebra from \eqref{E:noeth:gwa}. Prototypical examples of GWAs are (quantum) planes, (quantum) Weyl algebras, $U(\mathfrak{sl}_2)$, its quantum deformation $U_q(\mathfrak{sl}_2)$ and their infinite-dimensional primitive quotients, and Noetherian (generalized) down-up algebras. Amongst these, the plane (i.e., the commutative polynomial ring in two variables) is $\A_0$ and the Weyl algebra is $\A_1$.

\begin{theorem}[{\cite[Thm.\ 8.2]{BLO15tams}}]
	Assume  $h \not \in\FF$. Then the algebra $\A_h$ is not isomorphic to a generalized Weyl algebra over a polynomial ring in one variable. 
\end{theorem}

Turning to automorphisms, in \cite{BLO15tams} we give an exact description of the automorphism group of $\A_h$ over arbitrary fields $\FF$, describe the invariants in $\A_h$ under the automorphism group and investigate the analogue of Dixmier’s conjecture for $\A_h$ in case $\degg h \geq 1$.

For  $f\in\FF[x]\subseteq \A_h$ there is  $\phi_f\in\aut (\A_h)$  defined by 
\begin{equation}\label{E:autos_Ah:fromlnd} \phi_f (x)=x, \qquad \phi_f (y)=y+f(x).  \end{equation} 
Furthermore, $\{ \phi_f \mid f\in\FF[x]  \} \cong \left( \FF[x], +\right)$ is a subgroup of $\aut (\A_h)$. In case $h=1$ (Weyl algebra), $\sigma : \A_{1}\rightarrow\A_{1}$, $x \mapsto -y$, $y \mapsto x$ gives an automorphism of $\A_1$ of order $4$ and we have:

\begin{theorem}[Dixmier '68 ($\chara (\FF)=0$), Makar-Limanov '84 ($\chara (\FF)>0$)]
	$\aut(\A_{1})$ is generated by $\FF[x]$ and $\sigma$. 
\end{theorem}
Note that a similar result holds for $\aut (\A_0)$, by Jung ($\chara (\FF)=0$) and Van der Kulk ($\chara (\FF)>0$).

Assume now that $\degg h\geq 1$. Consider the automorphisms:
$$  \tau_{\alpha,\beta}(x) = \alpha x + \beta, \qquad \tau_{\alpha,\beta}(y) = \alpha^{\mathsf{deg}h-1} y. $$
for $\ds  (\alpha,\beta)\in \mathbb P = \{(\alpha,\beta) \in \FF^* \times \FF \mid  h(\alpha x +\beta) = \alpha^{\degg h} h(x)\} .$

The following are subgroups of $\aut(\A_h)$:
$$\tau_{\mathbb P} :=\{ \tau_{\alpha,\beta} \mid (\alpha,\beta) \in \mathbb P\} \quad \text{and} \quad  \tau_{1,\mathbb{G}} = \{ \tau_{1, \nu} \mid (1,\nu) \in \mathbb{P} \}.$$

\begin{theorem}[{\cite[Sec.\ 8]{BLO15tams}}]
	Assume that $\degg h\geq 1$ and $\FF$ is arbitrary.
	\begin{enumerate}[label=(\alph*)]
		\item $\FF[x]$ is a normal subgroup of $\aut (\A_h)$, and
		$\displaystyle \aut(\A_h) = \FF[x] \rtimes \tau_{\mathbb P}$.
		\item $\tau_{1,\mathbb{G}}$ is isomorphic to a finite subgroup of $(\FF,+)$, which  is trivial when $\chara(\FF) = 0$.
		\item $\FF[x] \rtimes \tau_{1,\mathbb G}$ is a normal subgroup of $\mathsf{Aut}_\FF(\A_h)$.
	\end{enumerate} 
\end{theorem}

The derivations of the algebras $\A_h$ were studied in detail over fields of arbitrary characteristic in \cite{BLO15ja}, where a very interesting connection with the Witt algebra was obtained. We will report on this in Section~\ref{S:hh} on $\hoch^\bullet(\A_h)$. However, we point out that the automorphisms \eqref{E:autos_Ah:fromlnd} arise in the form $\phi_f=e^{D_f}$, where $\pb{D_f\mid f\in\FF[x]}$ is an abelian Lie algebra of locally nilpotent derivations of $\A_h$. This motivated our study in \cite{KLM21} of the locally nilpotent derivations of $\A_h$ and more generally, in case $\chara(\FF)>0$, of Hasse--Schmidt higher derivations of $\A_h$. In fact, in \cite{KLM21} we connected the theorems of Rentschler \cite{rR68} and Dixmier \cite{jD68} on locally nilpotent derivations and automorphisms of the polynomial ring $\A_0$ and of the Weyl algebra $\A_1$ by establishing the same type of results for the family of algebras $\A_h$.

\subsection{Generalized down-up algebras}\label{SS:isoautder:gdua}

Motivated by the combinatorics of up and down operators on posets and in particular Stanley's $r$-differential posets (see Subsection~\ref{SS:intro:wa:comb}), Benkart and Roby defined down-up algebras in their seminal paper \cite{BR98}. 

\begin{definition}\label{D:def:dua}
	Given arbitrary fixed complex constants $\alpha, \beta, \gamma$, the down-up algebra $A(\alpha, \beta, \gamma)$ is the unital associative algebra generated by two generators $d$ and $u$, with relations:
	\begin{align*}
		d^2 u=\alpha dud+\beta ud^2+\gamma d \q\q \text{and}\q\q du^2=\alpha udu+\beta u^2d+ \gamma u. 
	\end{align*}
\end{definition}

Any down-up algebra with $(\alpha, \beta)\neq (0,0)$ is isomorphic to a Witten $7$-parameter deformation of $U(\mathfrak{sl}_{2})$ and they are also related to Le Bruyn's conformal $\mathfrak{sl}_{2}$ enveloping algebras~\cite{lLB95}.

Two of the most remarkable examples of down-up algebras are $A(2,-1,-2)=U(\mathfrak{sl}_{2})$ and 
$A(2,-1,0)=U(\mathfrak{h})$, the enveloping algebras of the $3$-dimensional complex simple Lie algebra 
$\mathfrak{sl}_{2}$ and of the $3$-dimensional nilpotent, non-abelian Heisenberg Lie algebra 
$\mathfrak{h}$, respectively. These algebras have a very rich structure and representation theory which has been extensively studied, having an unquestionable impact on the theory of semisimple and nilpotent Lie algebras. Nevertheless, a precise description of their symmetries, as given by the understanding of their automorphism group, is yet to be obtained (see~\cite{jD68, jD73} and~\cite{aJ76, jA86}).

Generalized down-up algebras were introduced by Cassidy and Shelton in~\cite{CS04} as a generalization of down-up algebras.

\begin{definition}\label{D:def:gdua}
	Let $\FF$ be an arbitrary field. The generalized down-up algebra $L(v,r,s,\gamma)$ is the unital associative $\FF$-algebra generated by $d$, $u$ and $h$ with defining relations
	\begin{equation*}
		dh-rhd +\gamma d =0, \quad \quad hu-ruh +\gamma u=0 \quad \mbox{and}\quad du-sud +v(h)=0,
	\end{equation*}
	where $r, s, \gamma\in\FF$ and $v\in\FF[h]$.
\end{definition}

Generalized down-up algebras include all down-up algebras, as long as the polynomial $h^2-\alpha h-\beta$ has roots in $\FF$. Indeed, if $\alpha=r+s$ and $\beta=-rs$, then it is easy to see that $A(\alpha, \beta, \gamma)\simeq L(-h,r,s,-\gamma)$. Conversely, any generalized down-up algebra $L(v,r,s,\gamma)$ with $\degg v=1$ is a down-up algebra. Generalized down-up algebras include also the algebras \emph{similar to the enveloping algebra of $\mathfrak{sl}_{2}$} defined by Smith~\cite{spS90}, Le Bruyn's \emph{conformal $\mathfrak{sl}_{2}$ enveloping algebras}~\cite{lLB95} and Rueda's algebras \emph{similar to the enveloping algebra of $\mathfrak{sl}_{2}$}~\cite{sR02}. 

By realizing $L(v,r,s,\gamma)$ as a generalized Weyl algebra, in \cite{CL09} we were able to determine the automorphism group of all conformal Noetherian generalized down-up algebras over an algebraically closed field of characteristic $0$ such that $r$ is not a root of unity. Specializing to the particular case of down-up algebras, where the description is simplest, we obtained the following characterization.

\begin{theorem}\label{T:adua:dua:main2}
	Let $A=A(\alpha, \beta, \gamma)$ be a down-up algebra, with $\alpha=r+s$ and $\beta=-rs$. Assume that $\beta\neq 0$ and 
	that one of $r$ or $s$ is not a root of unity. The group $\aut (A)$ of algebra automorphisms of $A$ is described bellow.
	\begin{enumerate}
		\item If $\gamma=0$ and $\beta=-1$ then $\aut (A)\simeq\left(\FF^{*} \right)^{2}\rtimes\ZZ/2\ZZ$;
		\item If $\gamma=0$ and $\beta\neq-1$ then $\aut (A)\simeq\left(\FF^{*} \right)^{2}$;
		\item If $\gamma\neq0$ and $\beta=-1$ then $\aut (A)\simeq\FF^{*}\rtimes\ZZ/2\ZZ$;
		\item If $\gamma\neq0$ and $\beta\neq-1$ then $\aut (A)\simeq\FF^{*}$.
	\end{enumerate}
	In all cases, the $2$-torus $\left(\FF^{*} \right)^{2}$ acts diagonally on the generators $d$ and $u$, $\mu\in\FF^{*}$ acts as multiplication by $\mu$ on $d$ and as multiplication by $\mu^{-1}$ on $u$, and the generator of the finite group $\ZZ/2\ZZ$ interchanges $d$ and $u$.
\end{theorem}

\subsection{The case of \texorpdfstring{$\hqfg$}{Hq(f,g)}}\label{SS:isoautder:hqfg}

Moving over to the quantum generalized Heisenberg algebras $\hqfg$, in \cite{LR22jaa} we tackled the isomorphism problem for this class and found that the isomorphism relation can be phrased in very concrete geometric terms, very much like in~\cite{BJ01}. 

\begin{proposition}[\cite{LR22jaa}]\label{Prop7} 
	Let $q\in\FF$ and $f,g \in\FF[h]$. The following define isomorphisms of qGHA.
	\begin{enumerate}[label=\Roman*.,ref=\Roman*]
		\item For all $\alpha\in\FF$, $\tau_{\alpha}:\hqfg\longrightarrow \hqfg[f(h-\alpha)+\alpha, g(h-\alpha)]$, defined on the canonical generators by $x\mapsto  x$, $y\mapsto  y$ and $h\mapsto  h-\alpha$.\label{P:7:I}
		
		\item For all $\lambda\in\FF^*$, $\sigma_{\lambda}:\hqfg\longrightarrow \hqfg[\lambda f(\lambda^{-1}h), g(\lambda^{-1}h)]$, defined on the canonical generators by $x\mapsto  x$, $y\mapsto  y$ and $h\mapsto  \lambda^{-1}h$.\label{P:7:II}
		
		\item For all $\lambda, \mu\in\FF^*$, $\rho_{\lambda,\mu}:\hqfg\longrightarrow \hqfg[f,\lambda\mu g]$, defined on the canonical generators by $x\mapsto  \lambda^{-1}x$, $y\mapsto  \mu^{-1}y$ and $h\mapsto h$.\label{P:7:III}
	\end{enumerate}
\end{proposition}

\begin{theorem}[\cite{LR22jaa}]\label{T:iso}
	Assume $q\neq0$ and $\degg f>1$. Then $\hqfg\simeq \hqqfg$ if and only if $q=q'$ and $(f',g')$ is obtained from $(f,g)$ via transformations of types \ref{P:7:I}, \ref{P:7:II}, \ref{P:7:III}. 
\end{theorem}

It follows in particular that, in case $q\neq 0$ and $\degg f>1$, the parameter $q$, as well as the integers $\degg f$ and $\degg g$, are invariant under isomorphism, showing that qGHA are indeed a vast generalization of generalized Heisenberg algebras and generalized down-up algebras.

Concerning derivations, in \cite{sL17} we investigated the locally finite and the locally nilpotent derivations of generalized Heisenberg algebras $\hh(f)=\hh_1(f, f-h)$ (see Subsection~\ref{SS:rep:qgha}), assuming that $\degg f>1$ (as $\degg f\leq1$ corresponds to a generalized down-up algebra). We also obtained the following result on gradings of $\hh(f)$.

\begin{theorem}[{\cite[Cor.\ 4.2]{sL17}}]\label{T:lfd:zgrad}
	Assume that $\degg f>1$. Then for any $\ZZ$-grading of $\hh(f)$, there is an integer $\ell\in\ZZ$ such that, relative to that grading, $x$ has degree $\ell$, $y$ has degree $-\ell$ and $h$ has degree $0$.
\end{theorem}

In terms of automorphism groups, in \cite[Thm.\ 5.1]{sL17} we described the automorphism group of $\hh(f)$ and showed that it is always abelian: it is isomorphic to $\CC^*\times \mathsf{C}$, where $\mathsf{C}$ is a finite cyclic group whose order divides $(\degg f) -1$. These results were extended to quantum generalized Heisenberg algebras $\hqfg$ over arbitrary fields in \cite{LR22jaa} and an interesting phenomenon was noticed. Although, as long as $\chara(\FF)\neq\degg f$, the automorphism group of a quantum generalized Heisenberg algebra $\hqfg$ with $q\neq 0$ and $\degg f>1$ is abelian and does not depend on the parameter $q$ (yet its isomorphism class does), if we allow $\chara(\FF)=\degg f$ then we can obtain non-abelian automorphism groups.

\section{Ring theoretical results}\label{S:rt}

\begin{flushright}{\it
		``One of the most fascinating developments in ring theory in recent
		years is the way in which large parts of algebraic geometry can now be stated entirely
		in terms of commutative Noetherian rings (Grothendieck~\cite{aG60}). In the other direction this has led
		to new tools for classifying and investigating these rings; furthermore, these applications are no longer confined to Noetherian rings, and although commutativity
		is assumed as a rule, one suspects that even this is not always essential.''} \\
	P. M. Cohn, \cite[Introduction]{pC73}
\end{flushright}

Throughout our work, we have proved several ring theoretical results pertaining to several classes of algebras. Since this is generically more technical, we omit several details and definitions, limiting this section to a brief summary.

\subsection{Noncommutative unique factorization}\label{SS:rt:fact}

Let $R$ be a commutative integral domain. Then the following conditions are equivalent:
\begin{itemize}
	\item Every nonzero element of $R$ which is not a unit is uniquely a product of irreducible elements up to the order of the product and up to taking associate irreducible elements.
	
	\item Every nonzero element of $R$ which is not a unit is a product of irreducible elements and the irreducible elements are prime.
	
	\item Every nonzero prime ideal of $R$ contains some prime element.
	
	\item $R$ is a UFD.
\end{itemize}
In case $R$ is also Noetherian, we can add another equivalent condition to the above list:
\begin{itemize}
	\item All nonzero minimal prime ideals of $R$ are principal.
\end{itemize}

As an example showing the geometric meaning of the above, the coordinate ring $R$ of an algebraic variety is a UFD if and only if every subvariety of codimension $1$ determines a principal ideal of $R$ (i.e., the subvariety is a complete intersection).

In \cite{aC84}, Chatters  introduced the notion of a noncommutative Noetherian unique factorization domain (Noetherian UFD, for short), based on the following result.

\begin{theorem}\label{T:rtr:UFD:conds}
	Let $R$ be a  Noetherian prime ring with at least one height-one prime ideal. Then the following conditions are equivalent for $R$:
	\begin{enumerate}
		\item Every height-one prime ideal of $R$ is principal and completely prime;
		\item $R$ is a domain and every non-zero element of $R$ can be written as $cp_1 p_2 \cdots p_n$, for some $c\in \mathcal{C}$ and a finite sequence  $p_1,  \ldots, p_n$ of prime elements of $R$, where  $\mathcal{C}$ is the set of elements of $R$ which are regular modulo all height-one prime ideals. 
	\end{enumerate} 
\end{theorem}

If $R$ is a commutative Noetherian domain then the elements of $\mathcal{C}$ must be units. An example in which the elements of $\mathcal{C}$ are not all units is $\A_x$, the enveloping algebra of the two-dimensional nonabelian Lie algebra, as $x\in\mathcal{C}$.

\begin{definition}
	A Noetherian UFD is a Noetherian domain with at least one height-one prime ideal and satisfying one of the equivalent conditions of Theorem~\ref{T:rtr:UFD:conds}. 
\end{definition}

Chatters' definition was later extended in \cite{CJ86}. To simplify the new definition, we assume the rings we deal with have finite  Gelfand-Kirillov dimension.

\begin{definition}
	A Noetherian UFR is a Noetherian prime ring all of whose height-one prime ideals are principal. 
\end{definition}

\begin{example}
	The following are examples of  Noetherian UFRs.
	\begin{itemize}
		\item $R[X]$, if $R$ is a Noetherian UFR;
		\item $M_n(R)$, if $R$ is a Noetherian UFR;
		\item some group rings of  polycyclic-by-finite groups over a commutative Noetherian UFD;
		\item certain Ore extensions of type $R[X; \sigma]$ or $R[X; \delta]$, if $R$ is a Noetherian UFR;
		\item the algebra $\mathcal{O}_q(M_{m, n})$ of quantum matrices, with $q$ not a root of unity;
		\item the quantum algebra $U^+_q(\mathfrak{g})$, with $q$ not a root of unity;
		\item the algebra of regular functions $\mathcal{O}(G)$, for $G$ a semisimple complex algebraic group which is connected and simply connected;
		\item  $\mathcal{O}_q(G)$ over $\mathbb{C}$ with $q$ transcendental, where $G$ is as above.
	\end{itemize}
\end{example}

Combining a result of Conze \cite{nC74} with one of Chatters \cite{aC84} we have the following.

\begin{theorem}
	If $\mathfrak{L}$ is a finite dimensional complex Lie algebra which is either solvable or semisimple, then its enveloping algebra $U(\mathfrak{L})$ is a Noetherian UFD. 
\end{theorem}

For more properties and results on noncommutative UFR and UFDs, see \cite{GY16}, where these notions were used to obtain the main results of that paper.

Given the similarities between down-up algebras and enveloping algebras, the above result inspired the paper \cite{LL13} where we analyzed the height-one prime ideals of Noetherian generalized down-up algebras $L(v,r,s,\gamma)$ and determined, in terms of the defining parameters, when they are noncommutative Noetherian UFDs or UFRs. 

Moreover, by considering cases in which the parameters $r$ and $s$ are roots of unity, we obtained some insight into the behavior of enveloping algebras over fields of finite characteristic (see~\cite{aB07} and references therein). Indeed, our analysis yielded the following result, which shows that, for generalized down-up algebras, the distinction between a Noetherian UFR and a Noetherian UFD depends only on the existence of torsion in the multiplicative subgroup of $\FF^{*}$ generated by $r$ and $s$.

\begin{theorem}[\cite{LL13}]\label{mainUFD}
	Let $L=L(v, r, s, \gamma)$ be a generalized down-up algebra with $rs \neq 0$. Then $L$ is a Noetherian UFD if and only if $L$ is a Noetherian UFR and $\langle r, s \rangle$ is torsionfree.
\end{theorem}

The above result reduces the study to the classification of those generalized down-up algebras which are Noetherian UFRs.

\begin{theorem}[\cite{LL13}]\label{mainUFR}
	Let $L=L(v, r, s, \gamma)$ be a generalized down-up algebra with $rs \neq 0$. Then $L$ is a Noetherian UFR  except if $v\neq 0$ and one of the following conditions is satisfied:
	\begin{enumerate}
		\item $v$ is not conformal, $r$ is not a root of unity, and there exists $\zeta\neq \gamma/(r-1)$ such that $v(\zeta)=0$;
		\item $v$ is conformal, $\langle r, s\rangle$ is a free abelian group of rank $2$, and there exists $\zeta\neq \gamma/(r-1)$ such that $v(\zeta)=0$;
		\item $\gamma\neq 0$, $r=1$, $s$ is not a root of unity, and $v\notin \FF$.
	\end{enumerate}
\end{theorem}

Moving to the algebras $\A_h$ form Definition~\ref{D:rep:ah}, we have the following classification.

\begin{theorem}[{\cite[Sec.\ 7]{BLO15tams}}]
	$\A_{h}$ is a Noetherian UFR. Moreover:
	\begin{enumerate}[label=(\alph*)]
		\item  If $\chara(\FF) =  0$, then $\A_{h}$ is a Noetherian UFD.
		\item If $\chara(\FF) = p > 0$, then $\A_{h}$ is not a Noetherian UFD.
	\end{enumerate}
\end{theorem}

\subsection{Double Ore extensions}\label{SS:rt:doe}

The notion of a double Ore extension was introduced in \cite{ZZ08} in order to obtain new examples of regular Artin--Schelter algebras of global dimension $4$ and using this construction in \cite{ZZ09} the authors obtain $26$ such families of algebras. Given certain similarities with iterated Ore extensions, in \cite{CLM11} we obtained necessary and sufficient conditions for double Ore extensions to be iterated Ore extensions. Moreover, we studied conditions for the double Ore extensions to be left (right) Noetherian, domains, prime, semiprime left (right) Noetherian, and semiprime left (right) Goldie, partially answering some questions in \cite{ZZ08}. Given the technical nature of the definitions involved, the interested reader is referred to \cite{CLM11} for more details.

\subsection{Properties of quantum generalized Heisenberg algebras}\label{SS:rt:qgha}

Recall Definition~\ref{D:qgha}. The quantum generalized Heisenberg algebra $\hqfg$ has a PBW-type basis of the form $\pb{x^i h^j y^k \mid i, j, k\in \ZZ_{\geq 0}}$, as proved in \cite[Lem.\ 2.2]{LR22ca}.
If $\degg f\leq 1$ then $\hqfg$ is isomorphic to a generalized down-up algebra and thus, by \cite[Cor.\ 2.4]{CS04}, $\operatorname{GKdim}\hqfg=3$. In view of this result and the PBW-type basis $\pb{x^i h^j y^k \mid i, j, k\in \ZZ_{\geq 0}}$, it would be natural to conjecture that qGHA generally have GK-dimension $3$. This turns out to be false and in fact we proved the following.

\begin{corollary}[{\cite[Cor.\ 2.7]{LR22jaa}}]\label{P:basic:GK}
	Let $\hqfg$ be a qGHA. Then $\operatorname{GKdim}\hqfg=3$ if $\degg f\leq 1$ and $\operatorname{GKdim}\hqfg=\infty$ if $\degg f>1$.
\end{corollary}

This unexpected result allows us to plainly identify generalized down-up algebras as the qGHA such that $\degg f\leq1$.

\begin{corollary}[{\cite[Cor.\ 2.8]{LR22jaa}}]\label{gduaimpdf1}
	The quantum generalized Heisenberg algebra $\hqfg$ is isomorphic to a generalized down-up algebra if and only if $\degg f\leq1$.
\end{corollary}

The other direction in which our ring-theoretical study of quantum generalized Heisenberg algebras has led us was the property of being Noetherian. Again taking inspiration in down-up algebras $A(\alpha, \beta, \gamma)$, which are Noetherian if and only if $\beta\neq 0$ (the main result in \cite{KMP99}), and generalized down-up algebras $L(v, r, s, \gamma)$, which are Noetherian if and only if $rs\neq 0$ (by \cite[Prop.\ 2.5, Prop.\ 2.6]{CS04}, in case $\chara(\FF)=0$ and $\barr\FF=\FF$), we investigated the conditions under which quantum generalized Heisenberg algebras are Noetherian. 

\begin{proposition}[\cite{LR22jaa}]\label{P:noetherian} Let $\FF$ be an arbitrary field.
	A qGHA $\hqfg$ is right (or left) Noetherian if and only if $\degg f=1$ and $q\neq0$.
\end{proposition}

It follows form the above cited sources that any generalized down-up algebra has the property that it is Noetherian if and only if it is a domain. It is natural to wonder whether this property still holds for a qGHA. We give a negative answer to the above question. 

\begin{proposition}[\cite{LR22ca}]\label{P:basic:domain}
	The qGHA $\hqfg$ is a domain if and only if $\degg f\geq 1$ and $q\neq 0$.
\end{proposition}

Thus, the fact that for generalized down-up algebras being Noetherian is equivalent to being a domain could be thought of as an anomaly detected by the broader context of qGHA.

\section{Combinatorics}\label{S:comb}

The combinatorics of the Weyl algebra alluded to in Subsection~\ref{SS:intro:wa:comb} is quite rich and has motivated a lot of research along with that of its quantum deformation introduced in Subsection~\ref{SS:rep:qpqwa}. In \cite{BLO13}, studying the representation theory of the algebras $\A_h$ (see Definition~\ref{D:rep:ah}) over fields of positive characteristic, the powers $\seq{h\frac{d}{dx}}^n$ come up in the action formulae.

\begin{example}\label{Exam:comb:E1} To simplify the notation, set $\partial_x=\frac{d}{dx}$. For $f\in\FF[x]$ we have:
	\begin{align*}
		\seq{h\partial_x}^1(f) &=  \cred{1} f^{(1)}h  \\
		\seq{h\partial_x}^2(f) &=  \cred{1}f^{(2)}h^2 +\cred{1} f^{(1)} h^{(1)} h  \\
		\seq{h\partial_x}^3(f) &=  \cred{1}f^{(3)} h^3 + \cred{3} f^{(2)} h^{(1)} h^2 + \cred{1}f^{(1)} h^{(2)} h^2 + \cred{1}f^{(1)}( h^{(1)})^2 h \\
		\seq{h\partial_x}^4(f) &= \cred{1}f^{(4)}h^4 + \cred{6} f^{(3)} h^{(1)}h^3 + \cred{4} f^{(2)} h^{(2)} h^3 + \cred{7} f^{(2)}(h^{(1)})^{2} h^2 
		\\
		& \qquad 
		+ \cred{1}f^{(1)} h^{(3)} h^3 + \cred{4} f^{(1)} h^{(2)}h^{(1)} h^2 + \cred{1}f^{(1)} (h^{(1)})^{3} h\\
	\end{align*}%
\end{example}

\subsection{Universal polynomials}\label{S:comb:universal}

Our initial problem in \cite[Sec.\ 8]{BLO13} was to study the expressions above and in particular the coefficients displayed in color. In \cite{BLR20} we extended this problem, considering a formal viewpoint, taking an arbitrary ring $A$ and a derivation $\partial$ of $A$. For $h\in A$, set $h^{[k]}=\partial^k(h)$, for all $k\geq 0$.

\begin{theorem}[\cite{BLR20}]\label{T:univpolys:univprop}
	For any $n\geq 0$, there is a unique polynomial $U_n$ in the mutually commuting variables $\pb{y_i}_{i\geq 0}$ and the noncommutative variable $t$ such that, for any ring $A$, derivation $\partial$ of $A$ and central element $h$ in $A$,
	\begin{equation*}
		\left(h\partial\right)^n 
		= \eval{U_n}{y_i=h^{[i]},\, t=\partial}=U_n(h, h^{[1]}, h^{[2]}, \ldots; \partial),
	\end{equation*}
	as endomorphisms of $A$.
\end{theorem}

Thus, from Example~\ref{Exam:comb:E1} above, we get 
\begin{align*}
	U_1 & =y_0  t\\
	U_2 & = y_0 y_1 t + y_0^2 t^2,\\
	U_3 & = y_0 y_1^2 t + y_0^2 y_2 t  + 3 y_0^2 y_1 t^2 + y_0^3 t^3, \\
	U_4&=y_0^4t^4+6y_0^3y_1t^3+4y_0^3y_2t^2+7y_0^2y_1^2t^2+y_0^3y_3t+4y_0^2y_1y_2t+y_0y_1^3t,
\end{align*}
so that $(h \partial)^n$ is obtained from $U_n$ by specializing $t$ to $\partial$ and $y_k$ to $\partial^k (h)$.

The polynomials $U_n$ can be described easily in terms of partitions and our main focus is on the coefficients of $U_n$. Given an integer partition $\lambda=(\lambda_1,\lambda_2,\ldots,\lambda_\ell)$, we define $y_{\lambda}=\prod_{i=1}^\ell y_{\lambda_i}$.

\begin{proposition}[\cite{BLR20}, compare {\cite[Sec.\ 8]{BLO13}}]\label{P:recrel}
	Assume $n\geq 1$. There exist positive integers $c^n_\lambda$, where $\lambda$ runs through the set of integer partitions of size $0\leq |\lambda|< n$, such that
	\begin{equation*}
		U_n=\sum_{k=1}^n\sum_{\lambda\vdash n-k} c^n_\lambda y_0^{n-\ell(\lambda)}y_{\lambda}t^{k}. 
	\end{equation*}
	Additionally, the coefficients $c^n_\lambda$ satisfy the recurrence relation
	\begin{equation*}
		c^1_\emptyset=1, \quad c^{n+1}_\lambda=c^{n}_\lambda + \sum_{i=0}^{n-1} (m_{i}(\lambda)+1) c^{n}_{\lambda[i+1]},
	\end{equation*}
	where:
	\begin{itemize}
		\item $c^{n}_{\lambda[i]}=0$ if  $m_i(\lambda)=0$;
		\item if $m_i(\lambda)>0$, $\lambda[i]$ is obtained from $\lambda$ by subtracting $1$ from a part of $\lambda$ of size $i$;
		\item $c^{n}_{\lambda}=0$ if $\lambda\vdash m$ with $m\geq n$.
	\end{itemize}
\end{proposition}

The coefficients occurring in $U_n$ form a family $c^n_\lambda$ of nonnegative integers indexed by partitions $\lambda$ and displaying very interesting combinatorial properties. Table~\ref{table:c} below lists the values $c^n_{\lambda}$ for $n\leq 5$ and hints at the connections with Stirling numbers of both kinds and Eulerian numbers.

\begin{table}[htbp]
	\[
	\begin{array}{c|*{12}{c}} 
		\lambda & {\emptyset} & {\Yboxdim{4pt} \yng(1)} & {\Yboxdim{4pt} \yng(2)} & {\Yboxdim{4pt} \yng(1,1)} & {\Yboxdim{4pt} \yng(3)} & {\Yboxdim{4pt} \yng(2,1)} & {\Yboxdim{4pt} \yng(1,1,1)} & {\Yboxdim{4pt} \yng(1,1,1,1)}& {\Yboxdim{4pt} \yng(2,1,1)}& {\Yboxdim{4pt} \yng(2,2)}& {\Yboxdim{4pt} \yng(3,1)}& {\Yboxdim{4pt} \yng(4)}\\[5pt]
		\hline
		c^1_{\lambda}&1 &&&&&&&&&&&\\[0.1cm]
		c^2_{\lambda}&1 &1  &&&&&&&&&&\\[0.1cm]
		c^3_{\lambda}&1 &3  &1 &1   &&&&&&&&\\[0.1cm]
		c^4_{\lambda}&1 &6  &4 &7   &1&4  &1   &&&&&\\[0.1cm]
		c^5_{\lambda}&1 &10 &10& 25 &5& 30& 15 &1&11&4&7&1
	\end{array}
	\]
	\caption{The coefficients $c^n_{\lambda}$ of the polynomials $U_n$. Each partition $\lambda$ is represented by its Young diagram. This table conceals the signless Stirling numbers of the first kind (as the sums $\sum_{\lambda \vdash n-k} c^n_{\lambda}$), the Stirling numbers of the second kind (as the coefficients $c^{n}_{(1^{n-k})}$ indexed by one-column shapes) and the Eulerian polynomials (whose coefficients are the sums $\sum_{\ell(\lambda)=k} c^n_{\lambda}$).}\label{table:c}
\end{table}

\subsection{Interpretations of the universal polynomials}\label{S:comb:comp}

It turns out that the polynomials $U_n$ (or specializations of these) had appeared already in \cite{BF12}, where the authors remarked that ``\textit{the symbolic problem [of calculating $U_n$] is indeed difficult}''. In \cite{BLR20} we found simple ways of computing the polynomials $U_n$ and many interesting interpretations of the coefficients $c^n_{\lambda}$. We also described $U_n$ in terms of other combinatorial structures, including increasing trees and partial subdiagonal maps. 

\begin{theorem}[Umbral formula from \cite{BLR20}]
	$U_n$ is obtained by applying 
	to 
	\begin{equation*}
		\prod_{i=0}^{n-1} \left( x_{i}+ \cdots + x_1 + x_0\right)
	\end{equation*}
	the $\ZZ$-linear map $\ZZ[x_0, x_1,\ldots,x_n]\longrightarrow R\langle t \rangle$ 
	\[
	x_n^{\alpha_n} x_{n-1}^{\alpha_{n-1}} \cdots x_1^{\alpha_1}x_0^k \ 
	\overset{\cred{\mbox{\scriptsize lowering}}}{\underset{\cred{\mbox{\scriptsize exponents}}}{\longmapsto}} \ y_{\alpha_n} y_{\alpha_{n-1}} \cdots y_{\alpha_1}t^k.
	\] 
\end{theorem}

An interesting application of these polynomials to formal solutions to differential equations is the following.

\begin{example}[Solution to the differential equation $X'(u)=Y(X(u))$]\hfill\\
	Consider the differential equation
	\[ X'(u)=Y(X(u)),\] 
	with initial condition $X(0)=0$.

	Write
	$\ds
	Y(u)=\sum_{i\geq 0} y_i \frac{u^i}{i!}$. Then there is a unique (formal) solution
	\begin{equation*}
		X(u)=\sum_{n\geq 1} x_n \frac{u^n}{n!} 
	\end{equation*}
	where
	\begin{equation*}
		x_n=\eval{U_{n-1}(y_0, \ldots, y_{n-2};t)}{t^k\leftarrow y_k}.
	\end{equation*}
\end{example}

Using an interpretation of $U_n$ in terms of subdiagonal partial maps from $[n]$ to $[n]$, we obtained a closed formula.

\begin{corollary}[\cite{BLR20}]\label{C:closed:form:d:is:1}
	Let $n\geq k \geq 1$ and $\lambda$ be a partition of $n-k$. Then 
	\begin{equation*}
		c^n_\lambda =\sum_{i_1, \ldots, i_{n-1}}
		\prod_{j=1}^{n-1}
		\genfrac{(}{)}{0pt}{}{j-i_{1}-\cdots -i_{j-1}}{i_j},
	\end{equation*} 
	where the sum is carried over all sequences of nonnegative integers whose nonzero terms are the parts of $\lambda$.
\end{corollary}

\subsection{Combinatorial interpretations of the coefficients \texorpdfstring{$c_{\lambda}^n$}{cλn}}\label{S:comb:special}

There are some interesting specializations of the polynomials $U_n$ and their coefficients $c_{\lambda}^n$. We mention a few below.

\begin{itemize}
	\item $\ds \sum_{\mu \vdash n-1} c_{\mu}^n  =  (n-1)! \q \text{and} \q
	\sum_{r=0}^{n-1} \sum_{\lambda \vdash r}  c_\lambda^n    =  n!$
	\item $\ds \sum_{\lambda\vdash n-k}c^n_\lambda=c(n, k)$, the signless Stirling number of the first kind.
	\item $\ds c^n_{1^{n-k}}=
	\genfrac{\{}{\}}{0pt}{}{n}{k}$, the Stirling number of the second kind.
	\item $\ds U_n(1, 1, 0, \ldots, 0;1)=B_n$, the Bell number.
	\item $\ds \sum_{\lambda\vdash n-k} c_{\lambda}^{n} \prod_i q_{\lambda_i}=
	\genfrac{\{}{\}}{0pt}{}{n}{k}_{q,1}$, where 
	$\genfrac{\{}{\}}{0pt}{}{n}{k}_{q,d}$ is the generalized Stirling number and \\ $q_k=\nolinebreak q(q-1)\cdots (q-k+1)$ is the falling factorial.
	\item $\ds \sum_{\lambda: \ell(\lambda)=k} c_{\lambda}^n =A(n, k+1)$, the Eulerian number, counting the number of permutations $\pi\in S_n$ with $k$ descents.
\end{itemize}

As a corollary, we obtained a formula for the Stirling numbers of the second kind which is generally much simpler than the usual one obtained using the inclusion-exclusion principle.

\begin{corollary} The following formula holds.
	\[ 
	\genfrac{\{}{\}}{0pt}{}{n}{k}=\sum_{1\leq a_1< \cdots <a_{n-k}\leq n-1} (1+a_1 -1)\cdots (1+a_{n-k} -(n-k)).\] 
\end{corollary}

\subsection{Modular behavior of the coefficients \texorpdfstring{$c_{\lambda}^n$}{cλn}}\label{S:comb:modp}

Let $p$ be a prime. It is well known that most Stirling numbers of both kinds $c(p,k)$ and 
$\genfrac{\{}{\}}{0pt}{}{p}{k}$ vanish modulo $p$. We showed that this property is shared by the coefficients $c^p_\lambda$ and, more generally, by $c^n_\lambda$, in case $n$ is a power of $p$.

\begin{theorem}[\cite{BLR20}]\label{T:modp}
	For any prime $p$, prime power $n=p^m$ and partition $\lambda$ with $|\lambda|\neq n-1$ and $|\lambda|$ not a multiple of $p$, we have $c^n_\lambda \equiv 0\ \modd p$. In particular, if $|\lambda|\neq 0, p-1$, then $c^p_\lambda \equiv 0\ \modd p$.
\end{theorem}

Recalling that
\[
c(n,k)=\sum_{\lambda \vdash n-k} c_{\lambda}^n, \qquad  
\genfrac{\{}{\}}{0pt}{}{n}{k}=c_{1^{n-k}}^n 
\quad \text{and}\quad
\genfrac{\{}{\}}{0pt}{}{n}{k}_{q,1}= \sum_{\lambda\vdash n-k} c_{\lambda}^{n}  (q)_{\lambda},
\]
it follows that for any prime $p$, $n=p^m$ and $1<k<n$ not a multiple of $p$, the Stirling numbers of both kinds 
$c(n,k)$ and 
$\genfrac{\{}{\}}{0pt}{}{n}{k}$ as well as the generalized Stirling numbers 
$\genfrac{\{}{\}}{0pt}{}{n}{k}_{q,1}$ are multiples of $p$.

\subsection{Generalization: the normally ordered form of \texorpdfstring{$\seq{h\partial^{d}}^n$}{(h∂d)n} and differential operator rings}\label{S:comb:gend}

All of the above was generalized to describe the normally ordered form of $\seq{h\partial^{d}}^n$. This includes the polynomials $U_{n,d}$, the umbral formula, the coefficients $c^{n, d}_\lambda$, their closed formula, modular behavior and formulas involving the generalized Stirling numbers, such as the one below:
\begin{equation*}
	\genfrac{\{}{\}}{0pt}{}{n}{k}_{q,d}= \sum_{\substack{\lambda\vdash nd-k\\ \ell(\lambda)\leq n-1}} c_{\lambda}^{n,d} (q)_{\lambda}.
\end{equation*}

We have the following general result.
\begin{theorem}\label{T:gen:d}
	Let $A$ be a ring, $h$ a central element of $A$ and $\partial$ a derivation of $A$. Then, for all $n, d\geq 1$ we have the following normal ordering identity in the formal differential operator ring $A[z; \partial]$ (skew polynomial ring of derivation type):
	\begin{align*}
		(h z^d)^n= \eval{U_{n, d}}{y_i=h^{[i]},\, t=z}
		=\sum_{k=d}^{nd}\sum_{\lambda\vdash nd-k} c^{n, d}_\lambda h^{n-\ell(\lambda)}h^{[\lambda]}z^{k},
	\end{align*}
	where $h^{[\lambda]}=h^{[\lambda_1]}\cdots h^{[\lambda_\ell]}$, for $\lambda=(\lambda_1, \ldots, \lambda_\ell)$.
\end{theorem}

\section{Hochschild cohomology and the Gerstenhaber bracket}\label{S:hh}

Recall the following result, already used in these notes.

\begin{lemma}[{\cite[Lem.\ 2.2]{BLO15tams}}]\label{lem:poly}  Let $\FF$ be a field of arbitrary characteristic. Assume that $\A = \FF[x][y,\sigma,\delta]$ is an Ore extension with $\sigma$ an automorphism of $\FF[x]$.   Then $\A$ is isomorphic
	to one of the following:
	\begin{itemize}
		\item[{\rm (a)}]   a  quantum plane;
		\item[{\rm (b)}]   a quantum Weyl algebra;
		\item[{\rm (c)}]   the unital associative algebra $\A_h$, with $h = h(x) \in \FF[x]$ (see Definition~\ref{D:rep:ah}).  
	\end{itemize}
\end{lemma}

Quantum planes and quantum Weyl algebras have been extensively studied, even from the homological point of view (see e.g.\ \cite{GG14}). On the contrary, we have not found any comprehensive study of the Hochschild (co)homology of the family of algebras $\A_h$, with $\FF$ and $h  \in \FF[x]$ arbitrary. This has been carried out mostly in \cite{BLO15ja} and \cite{LS21} and in this section we summarize these results.

\subsection{The Lie algebra \texorpdfstring{$\hoch^1(\A_h)$}{HH1(Ah)}}\label{SS:hh:hh1}

The results for $\hoch^1(\A_h)$ differ according to the characteristic of the base field and we treat them separately.

\subsubsection{The case \texorpdfstring{$\chara(\FF)=0$}{char(F)=0}}\label{SSS:hh1:char0}

The paper \cite{BLO15ja} is devoted to the study of the Lie algebra $\hoch^1(\A_h)$ over an arbitrary field $\FF$. We assume always that $h\neq 0$. Define the monic polynomial $\pi_{h}= \frac{h}{\gcd (h, h')}$.

\begin{theorem}[{\cite[Thm.\ 5.13]{BLO15ja}}]\label{thm:hochdec}
	Assume $\chara (\FF)=0$.  Then 
	\begin{equation*}
		\hoch^1(\A_{h})=\mathsf{Z}(\hoch^1 (\A_h))\oplus [\hoch^1 (\A_h), \hoch^1 (\A_h)], \ \ \hbox{\rm where}  \end{equation*}
	\begin{equation*}
		\mathsf{Z}(\hoch^1 (\A_h))=\bigg\{ D_{r\frac{h}{\pi_h}}\,\bigg | \,\degg r <\degg \pi_{h} \bigg\}  \, \hbox{\rm  and} \,
		\dimm \mathsf{Z}\left(\hoch^1 (\A_h)\right) = \degg \pi_{h}.
	\end{equation*}
	Above, $\mathsf{Z}(\hoch^1 (\A_h))$ stands for the center of the Lie algebra $\hoch^1 (\A_h)$ and the derivation $D_f$ is defined by $D_f(x)=0$, $D_f(y)=f$, for $f\in\FF[x]$.
\end{theorem}

An infinite-dimensional Lie algebra which plays an important role in the description of $\hoch^1(\A)$ is the \textit{Witt algebra}. A confusion with terminology may arise here, since the term Witt algebra has been used in the literature to mean two different things: the complex Witt algebra is the Lie algebra of derivations of the ring $\CC[z^{\pm 1}]$, with basis elements $w_n=z^{n+1}\frac{d}{dz}$, for $n\in\ZZ$;
while over a field $\KK$ of characteristic $p>0$, the Witt algebra is defined to be the Lie algebra of derivations of the ring
$\KK[z]/(z^p)$, spanned by $w_n$ for $-1\le n \le p-2$. 
Here we are considering a subalgebra of the first one (defined over the field $\FF$): 
\begin{equation}\label{E:def:Witt}
	\W = \spann_\FF\{w_i \mid i \geq -1\}, 
\end{equation}
equipped with the Lie bracket $\lb{w_m, w_n}=(n-m)w_{m+n}$, for $m, n\geq -1$. It is easy to check that if $\chara(\FF)=0$, then $\W$ is a simple Lie algebra (\textit{cf.}\ \cite[Lem.\ 5.19]{BLO15ja}). For the sake of simplicity and in accordance with the usage in \cite{BLO15ja}, we will abuse terminology and refer to the algebra $\W$ defined above as the Witt algebra. To make the distinction clear, we will call the Lie algebra of derivations of $\FF[z^{\pm 1}]$, with basis $\{w_i\}_{i \in\ZZ}$, the \textit{full Witt algebra}.

Let $\pr_1, \dots, \pr_t$ be the distinct monic prime factors of $h$ and, for $k \geq 0$, assume that $\pr_1, \ldots, \pr_k$ are the ones which occur with multiplicity $>1$. (When $k = 0$, no factor has multiplicity $\geq 2$.) Consider the ideal 
\begin{equation*} 
	\mathcal N=\spann_\FF\{\ad_{r a_n}\mid r \in  \pr_1 \cdots \pr_k\FF[x], \ n \geq 0\} 
\end{equation*} 
of $[\hoch^1(\A_{h}), \hoch^1(\A_{h})]$, where the $a_n$ are certain elements in the Weyl algebra $\A_1$ which normalize $\A_h$.

\begin{corollary}[{\cite[Cor.\ 5.22]{BLO15ja}}]\label{cor:structureL}
	Assume $\chara(\FF) = 0$ and let $h$ be as above. Then the following hold:
	\begin{itemize}
		\item[{\rm (i)}]  $\mathcal N$ is the unique maximal nilpotent ideal of $[\hoch(\A_{h}), \hoch(\A_{h})]$ and the quotient  Lie algebra
		$[\hoch(\A_{h}), \hoch(\A_{h})]/\mathcal N$ is the direct sum of $k$ simple Lie algebras  
		\begin{equation*}
			[\hoch(\A_{h}), \hoch(\A_{h})]/\mathcal N \  \cong \  \left( (\FF[x]/\FF[x] \pr_{1}) \otimes \mathsf{W} \right)\oplus \cdots\oplus \left( (\FF[x]/\FF[x] \pr_{k}) \otimes \mathsf{W} \right),
		\end{equation*} 
		where $\mathsf{W}$ is the Witt algebra. 
		\item[{\rm (ii)}]  If $\alpha_{i}\leq 2$ for all $1\leq i\leq t$, then $\mathcal N=0$. 
		\begin{enumerate}
			\item[{\rm (a)}] If  $\alpha_{i}=1$ for all $i$, then $[\hoch(\A_{h}), \hoch(\A_{h})] = 0$.   
			\item[{\rm (b)}] If some $\alpha_i = 2$,
			then $[\hoch(\A_{h}), \hoch(\A_{h})] $ is the direct sum of simple Lie algebras obtained from the Witt algebra $\W$ by extending the scalars. 
		\end{enumerate}
		\item[{\rm (iii)}] If there is $i$ such that $\alpha_{i}\geq 3$, then $\mathcal N\neq 0$,  and $[\hoch(\A_{h}), \hoch(\A_{h})]$ is neither nilpotent nor semisimple.
	\end{itemize}
\end{corollary} 
\newpage

We obtain the following special cases.
\begin{corollary}
	Assume that $\chara(\FF) = 0$. 
	\begin{itemize}
		\item[] ($\A_1$) $\der(\A_1) = \inder(\A_1)$, so $\hoch^1(\A_1) = (0)$; 
		
		\item[] ($\A_x$)  $\der(\A_x) = \FF D_1 \oplus \inder(\A_x)$, so $\hoch^1(\A_x)=\FF D_1$;
		
		\item[] ($\A_{x^2}$)   $\displaystyle \hoch^1(\A_{x^2})  =  \FF D_{x} \oplus \mathsf{W}$;
		
		\item[] ($\A_{x^n}$,  $n \geq 3$)   $\displaystyle \hoch^1(\A_{x^n})/\mathcal N  =  \FF D_{x^{n-1}} \oplus \mathsf{W}$, where $\mathsf{W}$ is  the Witt algebra. 
	\end{itemize}
\end{corollary}

\subsubsection{The case \texorpdfstring{$\chara(\FF)=p>0$}{char(F)=p>0}}\label{SSS:hh1:charp}

This case is more intricate and quite more technical. Although in \cite{BLO15ja} this case is completely studied, we just single out a few structural properties which are easy to state. See \cite[Sec.\ 6]{BLO15ja} for full details.

Let  $$\mathsf{Res}: \der(\A_h) \rightarrow \der(\mathsf{Z}(\A_h))$$ be the restriction map  and \[ \overbar{\mathsf{Res}}: \hoch^1(\A_h)  \rightarrow \der(\mathsf{Z}(\A_h))\] be the induced map.

\begin{theorem}[{\cite[Thm.\ 6.17]{BLO15ja}}]
	Assume that $\chara(\FF) = p>0$. Then the following hold.
	\begin{enumerate}[label=(\alph*)]
		\item $\im\,\mathsf{Res}= \im\,\mathsf{\overbar{Res}}$ is a free $\mathsf{Z}(\A_h)$-submodule of $\der(\mathsf{Z}(\A_h))$ of rank 2, and $\der(\mathsf{Z}(\A_h))$ is the Witt algebra in 2 variables;
		\item $\hoch^1(\A_h)$ is a  free $\mathsf{Z}(\A_h)$-module if and only if $\gcd (h, h')=1$. 
		In this case, $\mathsf{\overbar{Res}}$ is an isomorphism onto the image.
	\end{enumerate}
\end{theorem}

\subsection{The Gerstenhaber algebra structure of \texorpdfstring{$\hoch^\bullet(\A_h)$}{HH.Ah}}\label{SS:hh:gerst}

The aim of our paper \cite{LS21} was to describe the structure---given by the Gerstenhaber bracket---of the Hochschild cohomology spaces $\hoch^\bullet(\A_h)$  as Lie modules over $\hoch^1(\A_h)$.

In general, for an algebra $A$, Gerstenhaber introduced in \cite{mG64} a (graded) Lie algebra structure on the Hochschild cohomology $\hoch^{\bullet}(A)=\bigoplus_{n\geq 0}\hoch^n(A)$. This comes from what is now called a left-symmetric algebra (or pre-Lie algebra) structure on the Hochschild complex (bar resolution) of $A$. Although this structure does not depend on the particular resolution of $A$ chosen, Gerstenhaber's construction is very specific of the Hochschild complex and to date, in spite of may important results in this direction (see e.g.\ \cite{sW19,mSA17,yV19,NW16}), there is no known simple way of obtaining this additional structure from an arbitrary resolution.

The degree zero cohomology $\hoch^0(\A_h)$ has been computed in \cite[Section 5]{BLO15tams} and $\hoch^1(\A_h)$ was described in \cite{BLO15ja}, both over arbitrary fields. As a consequence of \cite[Sec.\ 3]{LS21}, we know that $\hoch^i(\A_h) =0$ for all $i\geq 3$. So our efforts go towards computing $\hoch^2(\A_h)$ and the Gerstenhaber structure of $\hoch^{\bullet}(\A_h)$. For the latter, we will make use of the method introduced in \cite{mSA17}.

\subsubsection{The case \texorpdfstring{$\chara(\FF)=0$}{char(F)=0}}\label{SSS:hh:gerst:char0}

As seen above, the Hochschild cohomology $\hoch^{\bullet}(\A_h)$ can be made into a Lie module for the Lie algebra $\hoch^{1}(\A_h)$ of outer derivations of $\A_h$, under the Gerstenhaber bracket. In our setting this is especially interesting in case $\chara(\FF)=0$ and $\gcd(h, h')\neq 1$ as then the description of $\hoch^{1}(\A_h)$ is related to the Witt algebra and, as we shall see, the $\hoch^{1}(\A_h)$-Lie module structure of $\hoch^{2}(\A_h)$ can be described in terms of the representation theory of the Witt algebra and also of the Virasoro algebra (see below).

\begin{theorem}[{\cite[Cor. 3.11, Rem. 3.13]{LS21}}]
	Assume  $\chara(\FF)=0$. There are isomorphisms
	\begin{equation*}
		\hoch^2(\A_h) \cong \A_h/\gcd(h, h')\A_h\cong \D[y],
	\end{equation*}
	where $\D=\left(\FF[x]/\gcd(h, h')\FF[x]\right)$.
	In particular, $\hoch^2(\A_h)=0$ if and only if $h$ is a separable polynomial; otherwise, $\hoch^2(\A_h)$ is infinite dimensional. 
\end{theorem}

In what follows in this subsection, we will assume the following conditions and notation.
\begin{itemize}
	\item $\chara(\FF)=0$;
	\item $\overline{\FF}=\FF$ (unnecessary but simplifies notation and statements);
	\item $h=\pr_1^{\alpha_1}\cdots\pr_k^{\alpha_k}\pr_{k+1}\cdots \pr_t$ is the decomposition of $h$ into irreducibles, with $k\leq t$ and $\alpha_1, \ldots, \alpha_k\geq 2$;
	\item Without loss of generality, $k\geq 1$, as otherwise $\hoch^2(\A_h)=0$;
	\item  $\W$ is the Witt algebra.
\end{itemize}
Then:
\begin{itemize}
	\item $\hoch^1(\A_h) = \mathsf{Z}(\hoch^1(\A_h))  \oplus [\hoch^1(\A_h),\hoch^1(\A_h)]$;
	\item $\mathcal{R}=
	\mathsf{rad}(\hoch^1(\A_h))=
	\mathsf{Z}(\hoch^1(\A_h))\oplus \mathcal N$ is the largest nilpotent ideal of $\hoch^1(\A_h)$.
\end{itemize}
Thus,
\begin{align*}
	\hoch^1(\A_h)/\mathcal{R}\simeq{\W\oplus\cdots\oplus\W} \q \text{($k$ copies of the Witt algebra)}
\end{align*}
and 
\begin{align*}
	\hoch^2(\A_h)\cong \D_1[y]\oplus\cdots\oplus\D_k[y],
\end{align*}
where $\D_i=\faktor{\FF[x]}{\pr_i^{\alpha_i -1}}$.

Our main results from \cite[Sec.\ 6]{LS21} can be summarized as follows.

\begin{enumerate}
	\item There is a filtration by $\hoch^1(\A_h)$-submodules
	\begin{equation*}
		\hoch^2(\A_h)=P_0\supsetneq P_1\supsetneq \cdots \supsetneq P_{m_h-1}\supsetneq P_{m_h}=0,
	\end{equation*}
	where $m_h=\max\pb{\alpha_1, \ldots, \alpha_k}-1$.
	\item $S_i=P_i/P_{i+1}$ is completely reducible
	\begin{align*}
		S_i= \bigoplus_{j:\alpha_j\geq i+2} {V}_{ij},
	\end{align*}
	with ${V}_{ij}$ simple.
	\item As a vector space, ${V}_{ij}=\FF[y]$.
	\item $\lb{\mathcal R, S_i}=0$ so $S_i$ is naturally a 
	$\hoch^1(\A_h)/\mathcal{R}=\W_1\oplus\cdots\oplus\W_k$-module, where $\W_j\cong\W$.
	\item $\lb{\W_{j'}, {V}_{ij}}=0$ for $j'\neq j$ and the action of $\W_j$ on ${V}_{ij}$ is given by:
	\begin{equation*}
		w_m . y^\ell =(\ell-(m+1)\mu_{ij})y^{m+\ell}, \qquad \mbox{for all $m\geq -1$ and $\ell\geq 0$,}
	\end{equation*}
	where $\mu_{ij}=\frac{\alpha_j-i}{\alpha_j-1}$.
	\item As $\hoch^1(\A_h)$-modules, ${V}_{ij}\cong {V}_{i'j'}$ if and only if $(i, j)=(i', j')$.
\end{enumerate}

It follows in particular that the composition length of $\hoch^2(\A_h)$ is $\degg\seq{\gcd(h, h')}$ and the composition factors are pairwise non-isomorphic as $\hoch^1(\A_h)$-modules (but not necessarily as $\W$-modules). 

\begin{theorem}
	Assume that $\chara(\FF)=0$.
	Then $\hoch^2(\A_h)$ is a semisimple $\hoch^1(\A_h)$-module if and only if $h$ is not divisible by the cube of any non-constant polynomial. 
\end{theorem}

It is well known that the Hochschild cohomology and its Gerstenhaber structure is a derived invariant, thus we obtain the following interesting corollary.

\begin{corollary}
	Assume that $\chara(\FF)=0$ and $\overline{\FF}=\FF$.
	Let $\lambda(h)$ denote the partition encoding the multiplicities of the irreducible factors of $h$. If $\lambda(h)$ and $\lambda(g)$ are different partitions, then $\A_{h}$ is not derived equivalent to $\A_{g}$. 
\end{corollary}

\begin{remark}\label{R:virasoro}
	A Lie algebra related to the Witt algebra and of utmost importance in theoretical Physics is the Virasoro algebra, denoted by $\mathsf{Vir}$. It has basis $\{w_i \mid i \in\ZZ\}\cup\{ c \}$ over $\FF$, with brackets
	\begin{equation*}
		[c, \mathsf{Vir}]=0 \quad\mbox{and}\quad  \lb{w_m, w_n}=(n-m)w_{m+n}+\delta_{m+n, 0}\frac{m^3-m}{12}c,
	\end{equation*}
	for all $m, n\in\ZZ$. The composition factors of $\hoch^2(\A_h)$ can actually be naturally embedded into irreducible modules for the Virasoro algebra. These are the so-called intermediate series modules and it is a result of Mathieu~\cite{oM92} that a Harish-Chandra module for $\mathsf{Vir}$ is either a highest weight module, a lowest weight module or an intermediate series module. 
\end{remark}

\subsubsection{The case \texorpdfstring{$\chara(\FF)=p>0$}{char(F)=p>0}}\label{SSS:hh:gerst:charp}

As before with $\hoch^1(\A_h)$, the Hochschild cohomology $\hoch^\bullet(\A_h)$ is a bit involved and quite technical over fields of positive characteristic, although an explicit description of $\hoch^2(\A_h)$ is given in \cite[Thm. 3.21]{LS21}. Here, we show only one of our results which is easy to state (but quite tricky to prove) and a few examples.

\begin{theorem}[{\cite[Thm. 3.24]{LS21}}]
	Assume  $\chara(\FF)=p>0$ and let $\mathsf{Z}(\A_h)$ denote the center of $\A_h$. Then $\hoch^2(\A_h)$ is a free $\mathsf{Z}(\A_h)$-module if and only if $\gcd(h, h')=1$. In this case, $\hoch^2(\A_h)$ has rank one over $\mathsf{Z}(\A_h)$ and, moreover, $\hoch^\bullet(\A_h)$ is a free $\mathsf{Z}(\A_h)$-module.
\end{theorem}

\begin{example}\hfill\\
	\begin{enumerate}
		\item $\A_1$ is the Weyl algebra and
		\begin{align*}
			\mathsf{HH^2}(\A_1) \cong \mathsf{Z}(\A_1) x^{p-1}y^{p-1},
		\end{align*}
		a rank-one module over $\mathsf{Z}(\A_1)=\FF[x^p, y^p]$. 
		\item $\A_x$ is the universal enveloping algebra of the two-dimensional non-abelian Lie algebra and
		\begin{equation*}
			\mathsf{HH^2}(\A_x) \cong  \mathsf{Z}(\A_x)x^{p}y^{p-1},
		\end{equation*}
		again a rank-one module over $\mathsf{Z}(\A_x)$.
		\item $\A_{x^2}$ is the Jordan plane.
		\begin{itemize}
			\item Case 1: $p=2$:
			\begin{equation*}
				\mathsf{HH^2}(\A_{x^2}) \cong \D \oplus \D x\oplus \D x^2y\oplus \mathsf{Z}(\A_{x^2})x^3y.
			\end{equation*}
			\item Case 2: $p=3$:
			\begin{equation*}
				\mathsf{HH^2}(\A_{x^2}) \cong \D \oplus \D x^{2}y \oplus \mathsf{Z}(\A_{x^2})x^{4}y^{2}.
			\end{equation*}
			\item Case 3: $p>3$:
			\begin{align*}
				\mathsf{HH^2}(\A_{x^2}) &\cong \bigoplus_{j=0}^{p-1} \D x^{2j}y^j \oplus \mathsf{Z}(\A_{x^2})x^{2p+1}y^{p-1}.
			\end{align*}
		\end{itemize}
		
		Above, $\mathsf{Z}(\A_{x^2})=\FF[x^p, x^{2p}y^p]$ and $\D=\mathsf{Z}(\A_{x^2})/x^p\mathsf{Z}(\A_{x^2})$.
	\end{enumerate}
	
\end{example}

\subsection{Final remarks: \texorpdfstring{$\A_h$}{Ah} as a deformation of \texorpdfstring{$\A_0=\FF[x,y]$}{A0=F[x,y]} and the twisted Calabi-Yau property}\label{SS:hh:cy}

In the commutative case, the concept of a Calabi-Yau algebra arose in connection with Calabi-Yau manifolds and mirror symmetry, the prototypical example being the coordinate ring of an affine Calabi-Yau variety. In \cite{vG07arXiv}, this notion was extended for noncommutative algebras.

\begin{definition}
	An algebra $A$ is a $\nu$-twisted Calabi-Yau algebra of dimension $d\geq 0$ if the following hold:
	\begin{itemize}
		\item $A$ is homologically smooth, i.e., it admits a finitely generated projective resolution of finite length, as a bimodule over itself;
		\item 
		$$
		\operatorname{Ext}_{A^e}(A,A^e)=
		\begin{cases}
			0,& \mbox{if $i\neq d$}\\
			A^\nu&\mbox{if $i = d$},
		\end{cases}
		$$
		where $A^e=A\otimes A^{\rm op}$, $\nu\in\aut(A)$ and $A^\nu$ is the bimodule with right $A$-action twisted by $\nu$.
	\end{itemize} 
	The automorphism $\nu$ is called the Nakayama automorphism of $A$. (It is unique up to inner automorphisms.) A Calabi-Yau algebra (in the sense of Ginzburg \cite{vG07arXiv}) is a twisted Calabi-Yau algebra whose Nakayama automorphism is an inner automorphism. 
\end{definition}

By \cite{vdB98}, if  $A$ is  a $\nu$-twisted Calabi-Yau algebra of dimension $d$ then there is a twisted Poincar\'e duality between homology and cohomology:
\begin{align}\label{E:TPD}
	\mathsf{HH^{d-i}}(A)\stackrel{\cong}{\longrightarrow} \mathsf{HH_i}(A,  A^\nu).
\end{align}

The following result is a consequence of \cite[Thm.\ 3.3, Rmk.\ 3.4, (2.10)]{LLW14}.

\begin{theorem}
	The algebra $\A_h$ is $\nu$-twisted Calabi-Yau with Nakayama automorphism satisfying $\nu(x)= x$ and $\nu(y)=y+h'$. Thus, the twisted Poincaré duality \eqref{E:TPD} holds for $\A_h$.
\end{theorem}

To finish this section, we remark that the construction in Subsection~\ref{SS:intro:wa:deform} realizing  the Weyl algebra $\A_1$ as a deformation of the commutative polynomial algebra $\A_0$ under the Weyl--Groenewold product, generalizes to show that, in case $\chara(\FF)=0$, all members of the family $\A_h$ can be seen as deformations of the commutative polynomial algebra $\A_0$. For that purpose, consider the derivations $\phi=\frac{d}{dy}$, $\psi=h(x)\frac{d}{dx}$ of $\A_0$ and define, for $a, b\in\A_0$,
\begin{align*}
	a\star b=\sum_{n\geq 0} \frac{\phi^n(a)\psi^n(b)}{n!}\hbar^n.
\end{align*}
This induces an associative product on $\A_0\formal{\hbar}$ with
\begin{eqnarray*}
	x\star x =x^2, && y\star y=y^2,\\
	y\star x=yx+h(x)\hbar, && x\star y=xy.
\end{eqnarray*}
So 
\begin{align*}
	y\star x-x\star y=h(x)\hbar
\end{align*}
and setting $\hbar=1$ we retrieve all members of the family $\A_h$ as deformations of $\A_0$.

\section{Nonassociative algebras: degenerations and geometric classification}\label{S:naa}

In this section, ``algebra'' will mean ``nonassociative algebra'', i.e.\ a vector space $A$ equipped with a bilinear product which need not be associative or unital. Another deviation from our previous setting is that we will tacitly assume that our algebras are finite dimensional.

The geometric deformation theory of algebras was initiated by Gabriel and others in the 1960s, inspired by the deformation theory of complex varieties in the work of Kodaira and Spencer \cite{KS58}. See also \cite{mG64}, \cite{NR66} and \cite{fF68}, where this was named ``algebraic geography''. The relationship between geometric features of the variety (such as irreducibility, dimension, smoothness) and the algebraic properties of its points brings novel geometric insight into the structure of the variety, its generic points and degenerations. 

Specifically, given a complex finite-dimensional vector space $V$ with a fixed basis $\pb{e_{1},\ldots ,e_{n}}$, the set of possible structure constants $\seq{c_{ij}^k}_{1\leq i,j,k\leq n}$ of an algebraic structure on $V$ (e.g.\ Lie, Jordan, Leibniz or anticommutative algebra), defined by
\begin{align*}
	e_i e_j=\sum_{k=1}^n c_{ij}^k e_k,
\end{align*}
forms an algebraic variety $X$ (under the Zariski topology). Algebraic properties of the points in $X$ often correspond to nice geometric features of $X$ and vice-versa. The general linear group acts on $X$ by coordinate change and the orbits form the isomorphism classes within the algebraic structure. Then, the orbit closures describe the degenerations in the variety: an algebra $B$ in the orbit closure of $A$ (but not isomorphic to $A$) is said to be a (proper) \textit{degeneration} of $A$ and we write $A\to B$. Moreover, an algebra whose orbit is open is said to be \textit{rigid} and its orbit closure is thus an irreducible component of the variety.

We name the classification of algebras by isomorphism in a given variety of algebras (of a given fixed dimension) defined by polynomial identities the ``algebraic classification''; it corresponds to the classification of orbits in the variety under the base-change action of the general linear group. Then the ``geometric classification'' is the study of the orbit closures, corresponding degenerations and rigid algebras and the determination of the irreducible components of the variety.

\begin{example}
	The following well-known varieties of $n$-dimensional algebras are defined by polynomial identities: 
	\begin{enumerate}
		\item Anticommutative algebras
		\begin{equation*}
			c_{ij}^k+c_{ji}^k=0, \quad 1\leq i, j, k\leq n.
		\end{equation*}
		\item Associative algebras
		\begin{equation*}
			\sum_{\ell=1}^n c_{ij}^\ell c_{\ell k}^m- c_{i\ell}^m c_{j k}^\ell=0, \quad 1\leq i, j, k, m\leq n.
		\end{equation*}
		
		\item Lie algebras
		\begin{align*}
			\sum_{\ell=1}^n c_{ij}^\ell c_{k\ell}^m+ c_{jk}^\ell c_{i\ell}^m+c_{ki}^\ell c_{j \ell}^m=0 \quad\mbox{and}\quad c_{ij}^k+c_{ji}^k=0, \quad 1\leq i, j, k, m\leq n.
		\end{align*}
	\end{enumerate}
\end{example}

\subsection{Central extensions and the algebraic classification}\label{SS:naa:ce}

The study of central extensions within a variety of algebras plays an important role in the classification problem in such varieties. Skjelbred and Sund \cite{ss78} devised a method for classifying nilpotent Lie algebras, making crucial use of central extensions, and it has since been adapted to many other varieties of algebras, including associative, Malcev, Jordan, Leibniz, and many others. For example, in Remark~\ref{R:virasoro} we have encountered the Virasoro Lie algebra, which is a central extension of the (full) Witt algebra.

Fix a variety $\mathcal M$ of algebras defined by polynomial identities, let $A$ be an algebra in $\mathcal M$ and let $V$ be a vector space. A $2$-cocycle of $A$ with values in $V$ is a bilinear map $\theta\colon { A}\times { A}\longrightarrow { V}$ satisfying the set of identities of $\mathcal M$ (see \cite[Sec.\ 1]{KLP20} for details). The set of $2$-cocycles as above is a vector space denoted by $\mathsf{Z}_{{\mathcal M}}({ A},{ V})$. Let $f\in\Hom_\FF(A,V)$. Then $\delta f=f\circ\mu\in\mathsf{Z}_{{\mathcal M}}({ A},{ V})$, where $\mu:{ A}\times { A}\longrightarrow { A}$ is the multiplication map. We define the space of $2$-coboundaries as $\mathsf{B}  ({ A},{ V})=\pb{\delta f \mid  f\in\Hom_\FF(A,V) }$. The {\it second cohomology} of $A$ with values in $V$ is the quotient space $\mathsf{H}^2_{\mathcal M}({ A},{ V})=\mathsf{Z}_{{\mathcal M}}({ A},{ V})/\mathsf{B}  ({ A},{ V})$.

For every bilinear map $\theta\colon { A} \times { A} \longrightarrow { V}$, we can define the algebra ${ A}_{\theta}={ A}\oplus { V}$ with product $[x+v,y+w]_{\theta}=xy+\theta(x,y)$. Then it follows (see \cite[Lem.\ 1.1]{KLP20}) that ${ A}_{\theta}$ is in the variety ${\mathcal M}$ if and only if $\theta\in\mathsf{Z}_{\mathcal M}(A, V)$.

The group $\aut(A)$ of algebra automorphisms of $A$ acts on the space $\mathsf{Z}_{\mathcal M}(A, V)$ by $\phi\cdot\theta(x,y)=\theta (\phi(x),\phi(y))$, for $\phi\in\aut(A)$ and $\theta\in\mathsf{Z}_{\mathcal M}(A, V)$. Moreover, it is easy to see that $\mathsf{B}  ({ A},{ V})$ is stabilized by this action so there is an induced action of $\aut(A)$ on $\mathsf{H}^2_{\mathcal M}({ A},{ V})$. The orbit of an element $\theta\in \mathsf{Z}_{\mathcal M}(A, V)$ will be denoted by $\operatorname{\mathsf{Orb}}(\theta)$, and similarly for $[\theta]\in\mathsf{H}^2_{\mathcal M}({ A},{ V})$.

Recall that the annihilator of an algebra $A$ is $\ann (A)=\pb{x\in{A}\mid xA + Ax=0}$.
For $\theta\in\mathsf{Z}_{\mathcal M}(A, V)$, define $\ann(\theta)=\pb{x\in A\mid \theta(x, A)+\theta ( A, x)=0}$. Next we remark that any algebra $A$ of the variety ${\mathcal M}$ with nonzero annihilator is isomorphic to a central extension of some other suitable algebra ${A'}$ in ${\mathcal M}$. 

\begin{lemma}[{\cite[Lem.\ 5]{HAC16}}]\label{lm:first}
	Let ${A}$ be an $n$-dimensional algebra in the variety ${\mathcal M}$ such that $\dim(\ann({ A}))=s\neq0$. Then there exist, up to isomorphism, a unique $(n-s)$-dimensional algebra ${A}'$ in $\mathcal M$ and a $2$-cocycle $\theta \in \mathsf{Z}_{\mathcal M}(A, V)$ for some vector space ${V}$ of dimension $s$, such that $\ann({A})\cap\ann(\theta)=0$,  ${A} \simeq {A'}_{\theta}$ and
	${A}/\ann({A})\simeq {A}'$.
\end{lemma}

Then, in order to decide when two algebras of the variety $\mathcal M$ with nonzero annihilator are isomorphic, it suffices to find criteria in terms of the cocycles, which is what is summarized below.

Fix a basis $\pb{e_{1},\ldots ,e_{s}}$ of ${V}$. Every cocycle $\theta\in \mathsf{Z}_{\mathcal M}(A, V)$ decomposes as $\theta=\sum_{i=1}^{s}\theta_i e_i$ with $\theta_i\in \mathsf{Z}_{\mathcal M}(A, \FF)$. Let $G_s(\mathsf{H}^2_{\mathcal M}({A},{\FF}))$ be the Grassmannian of $s$-dimensional linear subspaces of $\mathsf{H}^2_{\mathcal M}({A},{\FF})$ and set 
\begin{align*}
	E( {A},{V}) &=\pb{ {A}_{\theta }\mid \langle \left[ \theta _{1}\right],\dots,\left[ \theta _{s}\right] \rangle \in G_s(\mathsf{H}^2_{\mathcal M}({A},{\FF}))
		\q\text{and}\q 
		\bigcap\limits_{i=1}^{s}\ann(\theta _{i})\cap\ann({A}) =0}.
\end{align*}

The action of $\aut(A)$ on $\mathsf{H}^2_{\mathcal M}({A},{\FF})$ induces an action on $G_s(\mathsf{H}^2_{\mathcal M}({A},{\FF}))$ and we denote the orbit of $W\in G_s(\mathsf{H}^2_{\mathcal M}({A},{\FF}))$ by $\operatorname{\mathsf{Orb}}{(W)}$.

\begin{lemma}[{\cite[Lem.\ 17]{HAC16}}]
	Let ${A}_{\theta },{A}_{\vartheta }\in E( {A},{V}) $. 
	Then the algebras ${A}_{\theta }$ and ${A}_{\vartheta } $ of the variety ${\mathcal M}$ are isomorphic
	if and only if
	\[\operatorname{\mathsf{Orb}}(\langle \left[ \theta _{1}\right] ,
	\dots,\left[ \theta _{s}\right] \rangle) =
	\operatorname{\mathsf{Orb}} (\langle \left[ \vartheta _{1}\right] ,\dots,\left[ \vartheta _{s}\right]\rangle).\]
\end{lemma}

Suppose that the algebra $A$ is $n$-dimensional and consider the series
\[{A}^1={A}, \qquad \ {A}^{i+1}=\sum\limits_{k=1}^{i}{A}^k {A}^{i+1-k}, \qquad i\geq 1.\]
We say that ${A}$ is \emph{nilpotent} if ${A}^{i}=0$ for some $i \geq 1$. The smallest positive integer satisfying ${A}^{i}=0$ is called the  \emph{nilpotency index} of ${A}$.
Moreover, ${A}$ is called {\it null-filiform} if $\dim {A}^i=(n+ 1)-i$, for all $1\leq  i\leq n+1$.

By Lemma~\ref{lm:first}, any nilpotent algebra in a given variety can be obtained as a central extension of some nilpotent algebra of lower dimension in the same variety. The above results are a basis for the algebraic classification of such algebras. In particular, in \cite{KLP20} we classify central extensions of null-filiform associative algebras in several varieties of algebras containing the latter. This had been done in \cite{ACO17} for Leibniz algebras and in \cite{KLP20} we generalized it to  alternative, left alternative, Jordan, bicommutative, left commutative, assosymmetric, Novikov and left-symmetric  central extensions of null-filiform associative algebras.

In \cite{KKL20}, we completed the algebraic classification of all $6$-dimensional complex nilpotent anticommutative algebras. Concretely, we showed the following.

\begin{theorem}[{\cite[Thm.\ A]{KKL20}}]
	Up to isomorphism, the variety of $6$-dimensional complex nilpotent anticommutative algebras has infinitely many isomorphism classes, described explicitly in \cite[App.\ B]{KKL20} in terms of $14$ one-parameter families and $130$ additional isomorphism classes. 
\end{theorem}

In our paper \cite{KKL22}, we complete the algebraic classification of all complex $4$-dimensional nilpotent algebras, with no further restrictions. The final list, in \cite[Thm.\ 2]{KKL22},
has $234$ (parametric families of) isomorphism classes of algebras, $66$ of which are new in the literature.

\subsection{Geometric classification}\label{SS:naa:geom}

One of the main problems of the geometric classification of a variety of algebras is a description of its irreducible components. In~\cite{pG74}, Gabriel described the irreducible components of the variety of $4$-dimensional unital associative algebras and the variety of $5$-dimensional unital associative algebras was classified algebraically and geometrically by Mazzola~\cite{gM79}. Later, Cibils~\cite{cC87} considered rigid associative algebras with $2$-step nilpotent radical. Goze and Ancoch\'{e}a-Berm\'{u}dez proved that the varieties of $7$ and  $8$-dimensional nilpotent Lie algebras are reducible~\cite{GAB92}. The irreducible components of $2$-step nilpotent commutative associative algebras were described in~\cite{iS90}.

Often, the irreducible components of the variety are determined by the rigid algebras, although this is not always the case. Indeed Flanigan has shown in \cite{fF68} that the variety of $3$-dimensional nilpotent associative algebras has an irreducible component which does not contain any rigid algebras---it is instead defined by the closure of a union of a one-parameter family of algebras. 
We encounter similar situations in \cite{KKL20} and in \cite{KKL23}. 

Informally, although Theorem~\ref{T:B:20} shows that there is no single {\it generic} $6$-dimensional nilpotent anticommutative algebra, one can see the family ${\mathbb A}_{82}(\alpha)$ given below as the {\it generic family} in the variety. 

\begin{theorem}[{\cite[Thm.\ B]{KKL20}}]\label{T:B:20}
	The variety of $6$-dimensional complex nilpotent anticommutative algebras is irreducible of dimension $34$. It contains no rigid algebras and can be described as the closure of the union of $\mathrm{GL}_6(\mathbb{C})$-orbits of the following one-parameter family of algebras ($\alpha\in\mathbb{C}$):
	\begin{equation*}
		{\mathbb A}_{82}(\alpha) :\qquad   
		e_1e_2=e_3, \quad e_1e_3=e_4, \quad e_2e_5=\alpha e_6, \quad e_3e_4=e_5, \quad e_3e_5=e_6, \quad e_4e_5=e_6 .
	\end{equation*} 
\end{theorem}

In \cite{KKL23} we completely solve the geometric classification problem for nilpotent, commutative nilpotent and anticommutative nilpotent algebras of arbitrary dimension. We prove that the corresponding geometric varieties are irreducible, find their dimensions and describe explicit generic families of algebras which define each of these varieties. We also show some applications of these results in the study of the length of anticommutative algebras.

\begin{theorem}[{\cite[Thm.\ A]{KKL23}}]
	For any $n\ge 2$, the variety of all $n$-dimensional nilpotent algebras is irreducible and has dimension $\frac{n(n-1)(n+1)}{3}$.
\end{theorem}

Moreover, we show that the family $\mathcal{R}_n$ given in \cite[Def.\ 10]{KKL23} is generic in the variety of $n$-dimensional nilpotent algebras and inductively give an algorithmic procedure to obtain any $n$-dimensional nilpotent algebra as a degeneration from $\mathcal{R}_n$.

\begin{theorem}[{\cite[Thm.\ B]{KKL23}}]
	For any $n\ge 2$, the variety of all $n$-dimensional commutative nilpotent algebras is irreducible and has dimension $\frac{n(n-1)(n+4)}{6}$.
\end{theorem}

We show that the family $\mathcal{S}_n$ given in \cite[Def.\ 15]{KKL23} is generic in the variety of $n$-dimensional commutative nilpotent algebras and inductively give an algorithmic procedure to obtain any $n$-dimensional nilpotent commutative algebra as a degeneration from $\mathcal{S}_n$.

\begin{theorem}[{\cite[Thm.\ C]{KKL23}}]
	For any $n\ge 2$, the variety of all $n$-dimensional anticommutative nilpotent algebras is irreducible and has dimension $\frac{(n-2)(n^2+2n+3)}{6}$.
\end{theorem}

We show also that the family $\mathcal{T}_n$  given in \cite[Def.\ 33]{KKL23} is generic in the variety of $n$-dimensional anticommutative nilpotent algebras and inductively give an algorithmic procedure to obtain any $n$-dimensional nilpotent anticommutative algebra as a degeneration from $\mathcal{T}_n$.

The notion of length for nonassociative algebras has been recently introduced in~\cite{GK20}, generalizing the corresponding notion for associative algebras. Using the above result, we show in \cite[Cor.\ 39]{KKL23} that the length of an arbitrary (i.e.\ not necessarily nilpotent) $n$-dimensional anticommutative algebra is bounded above by the $n^{\text{th}}$ Fibonacci number, and prove that our bound is sharp.

\section*{Acknowledgment}

This work was Partially supported by CMUP, member of LASI, which is financed by national funds through FCT -- Funda\c c\~ao para a Ci\^encia e a Tecnologia, I.P., under the projects with reference UIDB/00144/2020 and UIDP/00144/2020.

This is an abridged version of our Habilitation thesis. We wish to thank the members of the Habilitation jury, namely Alberto Pinto, Nicolás Andruskiewitsch, Vladimir Bavula, Maria M. Clementino, Alberto Elduque, Pedro V. Silva and J. Toby Stafford, for the careful reading of a slightly extended version of these notes and for the kind and helpful questions and comments presented.

\def\cprime{$'$} \def\cprime{$'$} \def\cprime{$'$}

\EditInfo{August 1, 2023}{September 27, 2023}{Ana Cristina Moreira Freitas, Diogo Oliveira Silva, Ivan Kaygorodov, Carlos Florentino}

\end{document}